\documentclass{amsart}
\usepackage[english]{babel}
\usepackage{graphicx}
\usepackage[hidelinks]{hyperref}
\usepackage{amssymb}
\usepackage{amsmath, mathtools}
\usepackage[foot]{amsaddr}
\usepackage{amsfonts}
\usepackage{bbm}
\usepackage[dvipsnames, table]{xcolor}
\usepackage{pgfplots}
\usepackage{tikz}
\usepackage{tikzsymbols}
\usepackage{algorithm}
\usepackage{algpseudocode}

\usetikzlibrary{calc,trees,positioning,arrows,chains,shapes.geometric,%
    decorations.pathreplacing,decorations.pathmorphing,shapes,%
    matrix,shapes.symbols, decorations.markings, patterns,fit,quotes,%
    arrows.meta}
\usepackage[top=1in, bottom=1.25in, left=1.25in, right=1.25in]{geometry}
\usepackage{caption}
\usepackage{subcaption}
\usepackage{multirow}

\newtheorem*{theorem*}{Theorem}

\graphicspath{{images/}}

\input{Pack.sty}
\definecolor{maria}{HTML}{0090A0}

\begin{document}

\title[POD-\MakeLowercase{BASED} ROM\MakeLowercase{s for vb}OCP($\boldsymbol \mu$)\MakeLowercase{s}]{POD-based reduced order methods for optimal control problems governed by parametric partial differential equation with varying boundary control}
\author{Maria Strazzullo$^{\MakeLowercase{a},*}$, Fabio Vicini$^{\MakeLowercase{a},*}$}
\address{$^{\MakeLowercase{a}}$ INdAM research group GNCS member, Italy\newline $^{*}$ Politecnico of Torino, Department of Mathematical Sciences
``Giuseppe Luigi Lagrange", Corso Duca degli Abruzzi, 24, I-10129, Torino, Italy}

\begin{abstract}
In this work we propose tailored model order reduction for varying boundary optimal control problems governed by parametric partial differential equations. With varying boundary control, we mean that a specific parameter changes \emph{where} the boundary control acts on the system. This peculiar formulation might benefit from model order reduction. Indeed, fast and reliable simulations of this model can be of utmost usefulness in many applied fields, such as geophysics and energy engineering. However, varying boundary control features very complicated and diversified parametric behaviour for the state and adjoint variables. The state solution, for example, changing the boundary control parameter, might feature transport phenomena. Moreover, the problem loses its affine structure. It is well known that classical model order reduction techniques fail in this setting, both in accuracy and in efficiency. Thus, we propose reduced approaches inspired by the ones used when dealing with wave-like phenomena. Indeed, we compare standard proper orthogonal decomposition with two tailored strategies: geometric recasting and local proper orthogonal decomposition. Geometric recasting solves the optimization system in a reference domain simplifying the problem at hand avoiding hyper-reduction, while local proper orthogonal decomposition builds local bases to increase the accuracy of the reduced solution in very general settings (where geometric recasting is unfeasible). We compare the various approaches on two different numerical experiments based on geometries of increasing complexity.
\end{abstract}

\maketitle
\section{Introduction}
Parametric Optimal Control Problems (OCP($\bmu$)s) constrained to parametric partial differential equations (PDE($\bmu$)s) are widely used in connection with several engineering applications, such as geothermal and environmental analysis \cite{CARERE2021261, phd, GRENKIN20161243, Strazzullo1, Strazzullo3}, shape optimization \cite{delfour2011shapes, mohammadi2010applied} or the prediction of fluid flow and transport of contaminants \cite{BERRONE2018332, dede2007optimal, quarteroni2005numerical}. The main goal of OCP($\boldsymbol \mu$)s is to steer the expected behaviour of the PDE($\boldsymbol \mu$) towards a desired configuration by means of an external variable, called control. The control acts on the system in order to make the solution as close as possible to the desired state through a PDE($\boldsymbol \mu$)-constrained minimization problem.\
Each parameter $\bmu \in \mathcal P \subset \mathbb R^p$, represents a particular configuration, physical or geometrical, of the controlled system. For each new parameter, a different optimal solution is sought. In the aforementioned applications, OCP($\boldsymbol \mu$)s are usually related to time consuming activities where many simulations are required in a small amount of time. This task is not practicable in real-time contexts using standard discretization techniques that feature high computational costs.
For these reasons, in recent years, the research proposed different reduced order methods (ROMs) to reduce the total computational effort needed for simulating OCP($\bmu$)s, both at steady and time-dependent level \cite{dede2010reduced,Dede2012, hinze2005proper, Iapichino1,karcher2014certified, karcher2018certified,kunisch2008proper,Kunisch1999345,tesi,negri2013reduced,negri2015reduced,PichiStrazzullo,quarteroni2007reduced,Strazzullo2, StrazzulloRB,Strazzullo3, Troltzsch200983}.
ROMs strategies build a low-dimensional reduced space capable to describe the problem solutions varying with respect to $\boldsymbol \mu$ in a reliable and fast way. This approximation is based on the manipulation of some parametric instances of the \emph{high-fidelity} solution, i.e.\ solutions obtained by a standard discretization, in our case the Finite Element method (FE).

In this work we focus on steady Neumann boundary optimal control problems with a geometrical variation on the control action: the control variable plays a role on a portion of the boundary of the computational domain and this portion varies according to a specific geometric parameter $\mu_u \in \mathbb R$.
We are going to address this particular type of OCP($\boldsymbol \mu$) as varying boundary OCP($\boldsymbol \mu$) (vbOCP($\boldsymbol \mu$)).
Despite this parametric framework can be exploited to describe complex and interesting physical phenomena, to the best of our knowledge, this is the first time that this model is proposed.
As an example, the performances of heat exchangers can be improved by the variation of the baffle geometries, since the heat irradiated from baffles with different lengths has significant effects on the flow characteristics and heat transfer on the shell side \cite{BICER2020106417, YU2019351}.
The variation of the baffles seen as boundary controls on the heat exchanger, totally comply with the model at hand.
As a further example, the fluid flow inside porous fractured media is strongly influenced by the position and the size of the intersections among the fractures in the rock matrix, thus it is important to characterize the pressure field on each fracture changing the flow injected in the fracture intersections of different length \cite{BCPS,UQFR}.
Fractures moving and intersecting are a natural extension of vbOCP($\boldsymbol \mu$)s.

In literature, ROMs for standard geometric OCP($\boldsymbol \mu$) have been tackled by means of an affine transformation able to map the desired domain to a corresponding reference shape \cite{LassilaManzoniQuarteroniRozza2013a, rozza2013reduced, RLM}. However, the affinity hypothesis limits the applicability of the reduction in realistic geometrical configurations.
To overcome this problem, non-affine maps were proposed to introduce more complex deformations, \cite{JAGGLI2014263, LASSILA20101583}, coupled with tailored interpolation method to recover the linearity of the bilinear forms of the equations.
In this work, we are also interested in exploiting vbOCP($\boldsymbol \mu$) on complex domains and shapes, that can hardly be  recast in a reference domain. The main novelty of this contribution relies on:
\begin{itemize}
\item[\small{$\circ$}] the proposed vbOCP($\boldsymbol \mu$) model, which describes not a single control action but a \emph{varying control action} on the boundary of the considered domain;
\item[\small{$\circ$}] a first experimental analysis for ROMs approaches for vbOCP($\boldsymbol \mu$) based on tailored reduced strategies inspired from the ones used in nonlinear model order reduction, i.e. geometric recasting on a reference domain \cite{RozzaHuynhPatera2008} and local bases generation \cite{Amsallem2016,amsallem2012nonlinear,Borggaard20162155,Dihlmann2012156,drohmann2011adaptive, Eftang2011395,Eftang2011179,Haasdonk2011423, Maday2013A2417}.
\end{itemize}
Indeed, ROMs for vbOCP($\boldsymbol \mu$) turned out to be a very difficult task, that features complex parametric behaviour, wave-like-phenomena and non-affine structure. We believe that this work is a first step towards the applications of this model in real-life and interdisciplinary scenarios that naturally fit the proposed framework. We analyzed several reduction techniques and we compare them with standard proper orthogonal decomposition (POD) to propose some guidelines to deal with this peculiar optimal control problem.  
The comparisons are performed on two numerical tests of increasing geometrical complexity. The final goal is to provide a preliminary numerical investigation to analyze the best approach to be used in different geometrical settings of vbOCP($\boldsymbol \mu$)s.

The paper is outlined as follow. 
In Section~\ref{sec:problem} we present the vbOCP($\bmu$) continuous formulation and its high-fidelity formulation. Section~\ref{sec:rom} is devoted to introduce ROMs and the tailored reduction techniques needed to deal with vbOCP($\boldsymbol \mu$)s.
In Section~\ref{sec:results} we report the numerical tests in detail. We first test the approaches over a simple geometry and, then, we move our investigation towards more complex computational domains.
Conclusions are reported in Section~\ref{sec:conclusions}.
\section{Problem Formulation of vbOCP($\boldsymbol \mu$)s}
\label{sec:problem}
In this section we introduce vbOCP($\boldsymbol \mu$)s: i.e.\ varying boundary OCP($\bmu$)s, where the geometrical influence of the control, i.e.\ the portion of the boundary where the control acts, changes with respect to a parameter. We will focus on linear equations; however, the context can be easily generalized to other PDE($\bmu$)s. The discussion follows the formulation presented in other works on steady constrained optimization as \cite{karcher2014certified,karcher2018certified,negri2013reduced}. Indeed, even if in this problem the control action is different from standard geometrical OCP($\boldsymbol \mu$)s, it totally fits the standard framework. 

\subsection{Continuous formulation}
\label{sec:continuous}
\begin{figure}
    \centering
    \includegraphics[width=0.7\textwidth]{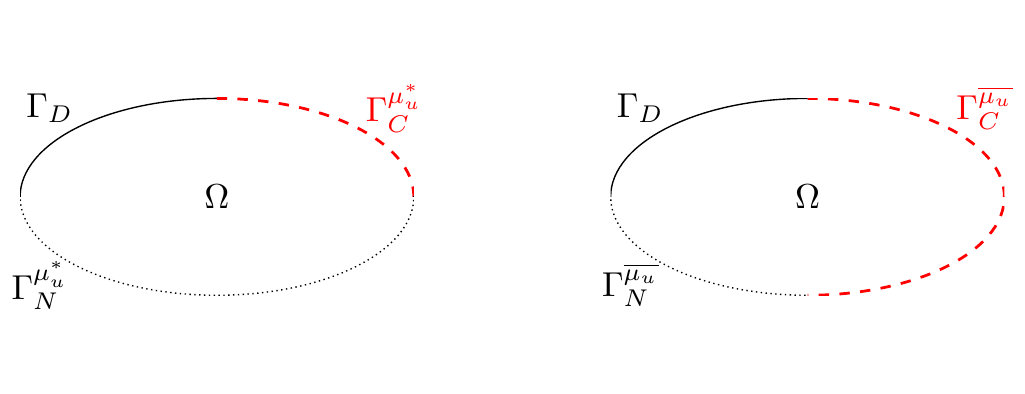}
\caption{Schematic representation of the computational domain and of the control action of a vbOCP($\boldsymbol \mu$) for different values of $\mu_u$, say $\mu_u = \mu_u^{*}$ (left) and $\mu_u = \overline{\mu_u}$ (right).}
\label{fig:domain_graph}
\end{figure}

Let us consider the parameter $\bmu \in \mathcal P \subset \mathbb R^p$, {$p \in \mathbb{N}$}: it represents physical and geometrical features of the problem at hand. In our specific setting, the parameter is of the form $(\mu^1, \dots, \mu^{p-1}, \mu_{u})$, where $\mu_u \in \mathbb R$ determines \emph{where} the boundary control is active.
We call $\Omega$ an open and bounded regular domain subset of $\mathbb R^2$. The portion of the boundary $\partial \Omega$ where Dirichlet conditions apply is $\Gamma_D$, while $\Gamma_N$ is the portion of $\partial \Omega$ where Neumann boundary conditions are considered. The Neumann boundary is split into $\Gamma_C^{\mu}$, the portion of the boundary where the control acts, and $\Gamma_N^{\mu}$, where homogeneous Neumann conditions are applied. 
The reason of the dependence on $\mu_u$ will be clarified in what follows. 
Moreover, let $y \in Y_g$ and $u \in U$ be the \emph{state} and the \emph{control} variables of the following spaces
$$
Y_g = \{y \in H^1(\Omega) \; | \; y = g \text{ on } \Gamma_D\} \text{ and } U = L^2(\Gamma_C^{\mu_u}),
$$
where $g \in H^{\frac{1}{2}}(\Gamma_D)$ is the value related to the Dirichlet boundary condition. From now on, we will assume $g \equiv 0$, without loss of generality: indeed, the non-homogeneous case can be dealt with a lifting procedure \cite{quarteroni2008numerical}. For the sake of notation, we call $Y = Y_0$.
In order to have a meaningful OCP($\bmu$), we still need to define, for $\Omega_{\text{obs}} \subset \Omega$, a \emph{desired state} $y_{\text{d}}(\bmu) \in L^2(\Omega_{\text{obs}})$. We are interested in the solution of:
\begin{equation*}
\min_{(y,u) \in Y \times U} 
\underbrace{\half \norm{y - y_{\text d}(\bmu)}_{L^2(\Omega_{\text{obs}})}^2 +\frac{ \alpha}{2} \norm{u}^2_{U}}_{J(y,u; \bmu)},
\end{equation*}
constrained to a PDE($\bmu$) of the following form
\begin{equation}
\label{eq:state_eq}
\begin{cases}
 {D}_a(\bmu)y = {f(\bmu)} & \text{in } {\Omega} , \\
\displaystyle y = 0 & \text{on  } {\Gamma_D}, \\
\displaystyle \dn{y} = u & \text{on  } {\Gamma_C^{\mu_u}}, \vspace{1mm}\\
 \displaystyle \dn{y} = 0 & \text{on  } {\Gamma_N^{\mu_u}}, \\
\end{cases}
\end{equation}
where $D_a(\bmu): Y \rightarrow Y\dual$ is a general differential operator, $f(\bmu) \in L^2(\Omega) \subset Y\dual$ is an external forcing term and $\alpha > 0$ is a penalization parameter for the control action. 
As already specified, the boundary \emph{where} the control acts changes with respect to the parameter $\mu_u$, as graphically represented in Figure~\ref{fig:domain_graph}. Namely, the parameter $\mu_u$ does not change the shape or the dimension of the domain by means of an affine transformation, but influences the measure of the curves $\Gamma_{C}^{\mu_u}$ and $\Gamma_{N}^{\mu_u}$. 

\begin{figure}
    \centering
    \includegraphics[width=0.7\textwidth]{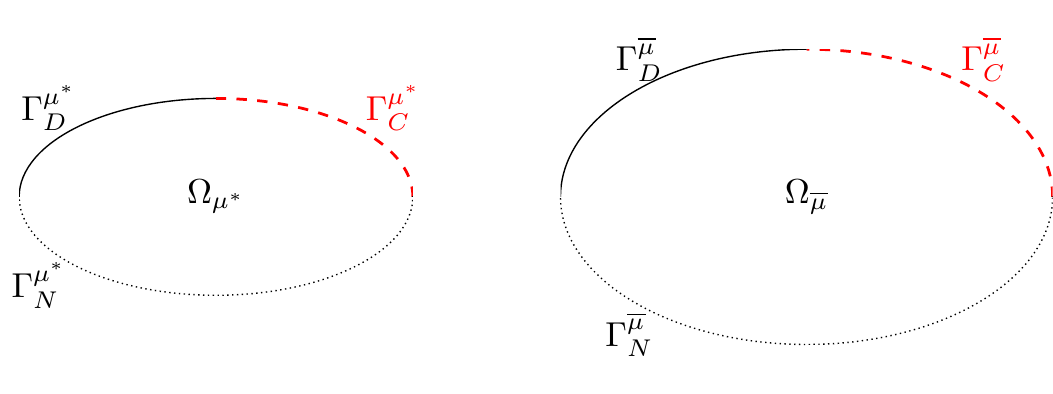}
\caption{Schematic representation of the computational domain and of the control action of a classical geometrical OCP($\boldsymbol \mu$) under affine transformation for different values of $\mu$, say $\mu = \mu^{*}$ (left) and $\mu = \overline{\mu}$ (right).}
\label{fig:domain_graph2}
\end{figure}

\emph{Parametric geometrical boundary control} was already investigated in several works, see for example \cite{negri2013reduced,Strazzullo2,StrazzulloRB,Zakia}, but in these cases, the geometrical parameter changes $\Omega \eqdot \Omega_{\mu}$ by means of an affine transformation. Consequently, the related portion of the boundary where the control acts varies, but not its relative measure, as depicted in Figure~\ref{fig:domain_graph2}. However, the herein problem is different: the parameter $\mu_{u}$ changes the nature of the problems itself, considering also limit scenarios, where $\Gamma_{N}^{\mu_u} \rightarrow \emptyset$ or $\Gamma_{C}^{\mu_u} \rightarrow \emptyset$. The setting we propose is denoted by the acronym vbOCP($\boldsymbol \mu$), as already specified. 

The consequences of this framework on the parametric reduction of the model will be widely discussed in the following sections. 
Let us focus on the problem formulation and the high-fidelity approximation of it.
The weak formulation of problem \eqref{eq:state_eq} reads: given $\bmu \in \Cal P$ find the pair $(y,u) \in Y \times U$ which verifies
\begin{equation}
\label{eq:state_eq_weak}
{a (y, q; \bmu)} =
{c(u ,q; \bmu)}
+ { \la G(\bmu), q \ra_{Y\dual, Y}} \quad \forall q \in Y,
\end{equation}
where the differential operator ${D}_a(\bmu)$ is represented by the bilinear form  $a \goesto{Y}{Y}{\mathbb R}$,  forcing terms are included in $ G(\bmu) \in Y\dual $ and, $c:U \times Y \rightarrow \mathbb R$ is defined by:
$$
{c(u ,q; \bmu)} \eqdot \int_{\Gamma_{C}^{\mu_u}}u q \; \text{ds},
$$
with $q \in Y$ meant in the sense of traces. \\
For the sake of notation, we drop the parameter dependence from the variables. Indeed, the value of $\bmu \in \mathcal P$ affects the optimal solution: namely, $y \eqdot y(\bmu)$ and $u \eqdot u(\bmu)$. We do the same for the desired state, i.e.\ $y_{\text{d}} \eqdot y_d(\bmu)$.


To state the well-posedness of the problem at hand, we assume that:
\begin{enumerate}
\item[(a)] \label{coercivity_a} $a(\cdot, \cdot; \bmu)$ is continuous and coercive;
\item[(b)] \label{continuity_of_c} $c(\cdot, \cdot; \bmu)$ is continuous.
\end{enumerate}
Thanks to (a) and (b), given a parameter $\bmu \in \mathcal P$ and a control $u \in U$, the weak \emph{state equation} \eqref{eq:state_eq_weak} is well-posed \cite{quarteroni2008numerical}. 
Moreover, combining (a) and (b)  with the definition of $Y$ and $U$ and $\alpha > 0$, for every parametric instance, we can state the well-posedness of the following \emph{weak minimization problem}: 
\begin{equation}
\label{eq:min_weak}
\min\{J(y, u; \bmu)  :  (y, u) \in Y \times U \text{ and \eqref{eq:state_eq_weak} is verified}\}.
\end{equation}
To solve  \eqref{eq:min_weak}, we rely on a Lagrangian approach, exploitable since \cite[Corollary 1.3]{hinze2008optimization} holds under our assumptions.
Thus, let us define the \emph{adjoint variable} $p \eqdot p(\bmu) \in Y$ and the Lagrangian functional as
\begin{align*}
 \Lg (y,u,p; \bmu) = J(y, u; \bmu)
+ \displaystyle a (y, p; \bmu)- c(u,p; \bmu)
-  \la G(\bmu), p \ra_{Y\dual, Y}.
\end{align*}
The minimizing pair $(y,u)$ is given by the solution of the following three-equations system:
 given $\bmu \in \Cal P$, find $(y, u, p) \in Y \times  U \times Y$ such that
\begin{equation}
\label{eq:optimality_system}
\begin{cases}
D_y\Lg(y, u, p; \bmu)[z] = 0 & \forall z \in Y,\\
D_u\Lg(y, u, p; \bmu)[v] = 0 & \forall v \in U,\\
D_p\Lg(y, u, p; \bmu) [q]= 0 & \forall q \in Y,\\
\end{cases}
\end{equation}
which in strong form reads: 
for a given $\bmu \in \mathcal P$ find $(y,u,p) \in Y \times U \times Y$ such that
\begin{equation}
\label{eq:strong_optimality_system}
\begin{cases}
\displaystyle y \chi_{\Omega_{\text{obs}}} + D_a(\bmu)\dual p =
y_d\chi_{\Omega_{\text{obs}}} & \text{ in } \Omega, \\
\alpha u -  p = 0 & \text{ in } \Gamma_C^{\mu_u}, \\
D_a (\bmu)y  = f (\bmu)& \text{ in } \Omega, \\
\displaystyle \frac{\partial y}{\partial n} = u & \text{ on }\Gamma_C^{\mu_u}, \vspace{1mm} \\ 
\displaystyle \dn{y} = \dn{p} = 0 & \text{ on } \Gamma_N^{\mu_u}, \\
y = p =0 &  \text{ on } \Gamma_D, \\
\end{cases}
\end{equation}
where we define $D_a(\bmu)\dual$ as the dual operator of $D_a(\bmu)$ and $\chi_{\Omega_{\text{obs}}}$ is the characteristic function of $\Omega_{\text{obs}}$. The first equation is called \emph{adjoint equation}, the second one \emph{optimality equation} and the third one is the state equation itself. It is clear that, thanks to the equality 
$\alpha u = p$ on $\Gamma_C^{\mu_u}$ (also in this case, we mean $p$ in the sense of traces), system \eqref{eq:strong_optimality_system} can be recast in: 
given $\bmu \in \Cal P$, find the pair $(y,p) \in Y \times Y$ such that
\begin{equation}
\label{eq:no_u}
\begin{cases}
\displaystyle y \chi_{\Omega_{\text{obs}}} + D_a(\bmu)\dual p =
y_d\chi_{\Omega_{\text{obs}}} & \text{ in } \Omega, \\

D_a (\bmu)y  = f (\bmu)& \text{ in } \Omega, \\
\displaystyle \frac{\partial y}{\partial n} = \frac{1}{\alpha}p & \text{ on }\Gamma_C^{\mu_u}, \vspace{1mm} \\ 
\displaystyle \dn{y} = \dn{p} = 0 & \text{ on } \Gamma_N^{\mu_u}, \\
y = p =0 &  \text{ on } \Gamma_D. \\
\end{cases}
\end{equation}
Namely, we consider only two variables and the control is recovered through a postprocessing procedure thanks to the condition $\alpha u -p = 0$ on the control boundary. 
System \eqref{eq:no_u} can be interpreted in a weak form as follows: given $\bmu \in \mathcal P$, find $(y,p) \in Y\times Y$ such that
\begin{equation}
\label{eq:weak_no_u}
\begin{cases}
m(y - y_{\text{d}}, w; \bmu) + a(w, p; \bmu) = 0 & \forall w \in Y,\\
\displaystyle a(y, q; \bmu) - \frac{1}{\alpha}c(p, q; \bmu) =  \la G(\bmu), q \ra_{Y\dual, Y}& \forall q \in Y,
\end{cases}
\end{equation}
where $m\goesto{Y_{\text{obs}}}{Y_{\text{obs}}}{\mathbb R}$ is the $L^2-$inner product restricted to the observation domain. \\
In the next section, we focus on the discretized version of the problem at hand.


\subsection{High-fidelity formulation}
\label{sec:hf}
To deal with system \eqref{eq:weak_no_u}, as already specified in the introduction, we use FE as high-fidelity  approximation.
Let us define a conforming and regular triangulation $\mathcal T$ on $\Omega$. Let
$Y\disc =
 Y \cap\Cal X_{r}$, with
$$
 \mathcal X_r = \{ y \in C^0(\overline \Omega) \; : \; y |_{K} \in \mbb P^r, \; \; \forall K \in \Cal T  \},
$$
where $K$ is an element of $\Cal T$ and $\mathbb P^r$ is the space of polynomials of degree at most $r$. \\
Any $y\disc$ and $p\disc$ in $Y\disc$ has a FE expansion of the form
$$
y\disc = \sum_{1}\disc y^{i} \phi_{i}, \hspace{1cm}
p\disc = \sum_{1}\disc p^{i} \phi_{i}, 
$$
where $\{\phi_i\}_{i=1}\disc$ is a basis for $Y\disc$ and $ y^{i}$, $p^i$ are (unknown) real numbers for $i=1,\dots, N_h$ (also in this case, we are omitting the $\boldsymbol \mu$-dependence of the FE coefficients).
We now perform a Galerkin projection in the FE spaces, solving
\begin{equation}
\label{eq:weak_no_u_FE}
\begin{cases}
m(y\disc - y_{\text{d}}, w; \bmu) + a(w, p\disc; \bmu) = 0 & \forall w \in Y\disc,\\
\displaystyle a(y\disc, q; \bmu) - \frac{1}{\alpha}c(p\disc, q; \bmu) =  \la G(\bmu), q \ra_{Y\dual, Y}& \forall q \in Y\disc.
\end{cases}
\end{equation}
Algebraically, we consider $\mathsf D_a(\bmu) = a(\phi_j, \phi_i; \bmu)$, for $i,j = 1, \dots, N_h$, as the stiffness matrix. Moreover,
 $\mathsf M_{{\text{o}}}(\bmu)$ is the mass matrix related to the $L^2$-inner product over $\Omega_{\text{obs}}$ and $\mathsf {C}(\bmu)_{ij} \eqdot c(\phi_j, \phi_i; \bmu)$, for $i,j = 1, \dots, N_h$ is the control matrix that acts over $\Gamma_C^{\mu_u}$. 
 We denote by $\mathsf y \in \mathbb R^{N_h}$ and $\mathsf p \in \mathbb R^{N_h}$ the unknown FE vectors for state and adjoint variables respectively. The same notation is used for the forcing term and the desired state, i.e. $\mathsf f = [f^1, \dots, f^{N_h}]^T$ and $\mathsf y_{\text{d}} = [y_{\text{d}}^1, \dots, y_{\text{d}}^{N_h}]$ with 
$$
f^i \eqdot \int_{\Omega} f(\bmu)\phi_i\; \text{dx} \quad \text{and} \quad y_{\text{d}}^i \eqdot  \int_{\Omega_{\text{obs}}} y_{\text{d}}(\bmu)\phi_i\; \text{dx},\ \forall i \in \{1, \dots, N_h\}.
$$
The algebraic formulation of the state equation reads 
\begin{equation}
\label{eq:alg_state_eq}
 \mathsf D_a(\bmu) \mathsf y- \frac{1}{\alpha} \mathsf C(\bmu) \mathsf p =  \mathsf f,
\end{equation}
while the adjoint equation has the following form:
\begin{equation}
\label{eq:alg_adj_eq}
 \mathsf M_{\text{o}}(\bmu) \mathsf y  + \mathsf D_a(\bmu)^T \mathsf p = \mathsf y_{\text{d}}. 
\end{equation}
Namely, combining equations \eqref{eq:alg_state_eq} and \eqref{eq:alg_adj_eq}, finding the minimizing FE solution $(y\disc, p\disc) \in Y\disc \times Y\disc$ translates in solving the saddle point system
\begin{equation}
\label{eq:one_shot}
\begin{bmatrix}
 \mathsf M_{\text{o}}(\bmu) &  \mathsf D_a(\bmu)^T  \\
\displaystyle \mathsf D_a(\bmu) & -\frac{1}{\alpha} \mathsf C(\bmu) \\
\end{bmatrix}
\begin{bmatrix}
\mathsf y \\
\mathsf p \\
\end{bmatrix}
=
\begin{bmatrix}
\mathsf y_{\text{d}} \\
\mathsf f \\
\end{bmatrix}.
\end{equation}
The two-equations system preserves the classical saddle point framework of optimal control problems: indeed,  
it is well-known that PDE($\bmu$) constrained optimization leads to such a peculiar structure both at the steady level \cite{bader2015certified,dede2010reduced,gerner2012certified,karcher2014certified,karcher2018certified,kunisch2008proper,negri2015reduced,negri2013reduced} and at the time-dependent one \cite{HinzeStokes,HinzeNS,hinze2008optimization,Stoll1,Stoll,Strazzullo2,Strazzullo3,StrazzulloRB}.

For the well-posedness of the presented saddle point structure, we refer the reader to \cite{Benzi}. We remark that the existence and uniqueness of an optimal solution can be proved by means of Brezzi Theorem \cite{Babuska1971,Brezzi74} as done in \cite{negri2013reduced} for the three-equations system \eqref{eq:strong_optimality_system}.
Another proving strategy exploits the Ne\v cas-Babu\v ska theory \cite{necas} as done in \cite{StrazzulloRB}. Both approaches consider the state and the adjoint solutions in $Y^{N_h}$. This is the reason why we assumed, from the very beginning, to work with the same space for both variables.

We are now interested in studying vbOCP($\boldsymbol \mu$) for several parametric instances: the high-fidelity system is usually not suited to analyze many parameters in a small amount of time. The FE problem suffers the high dimensionality of the optimality system, most of all when the physical phenomenon studied needs fine mesh resolutions. Thus, in the next Sections, we introduce ROMs for vbOCP($\boldsymbol \mu$).

\section{Reduced Order Methods for vbOCP($\boldsymbol \mu$)s}
\label{sec:rom}
We propose ROMs as a way to reliably solve vbOCP($\bmu$)s in a faster way. We introduce the general ideas of this modelling strategy and then we will move towards the algorithms we use to tackle the experiments presented in Section~\ref{sec:results}.
We describe three different approaches to deal with vbOCP($\bmu$)s, namely:
\begin{itemize}
\item[\small{$\circ$}] (\emph{Section~\ref{sec:POD}}) Proper Orthogonal Decomposition (POD).
\item[\small{$\circ$}] (\emph{Section~\ref{sec:LPOD}}) Local Proper Orthogonal Decomposition (L-POD).
\item[\small{$\circ$}] (\emph{Section~\ref{sec:geo-r}}) Geometric Recasting (Geo-R).
\end{itemize} 
The motivations behind the use of L-POD and Geo-R are discussed in Section~\ref{sec:transport}.
\subsection{Reduced Problem Formulation}
\label{off_on}
In this section we describe the basic ideas behind ROM strategies. To make concepts clearer, we explicit the parameter dependence of the variables.
Let us consider
\begin{equation}
    \label{eq:solution_manifold}
    \mathbb M = \{ (y(\boldsymbol{\mu}), p(\boldsymbol{\mu}))\;| \; \boldsymbol{\mu} \in \Cal P\},
\end{equation}
as the \emph{solution manifold}, i.e.\ the set of the optimal solutions varying with respect to the parametric instance. We here assume that $\mathbb M$ is smooth with respect to $\boldsymbol \mu \in \mathcal P$. Once defined {a high-fidelity} discretization, a \emph{discrete solution manifold} is defined analogously as
$$
    \mathbb M \disc= \{ (y\disc(\boldsymbol{\mu}), p\disc(\boldsymbol{\mu}))\;| \; \boldsymbol{\mu} \in \Cal P\}.
$$
Namely, $\mathbb M\disc\! \approx\! \mathbb M$ is a reliable surrogate of $\mathbb M$ when the high-fidelity approximation is a good representation of the the continuous model, i.e.\ when the mesh is fine enough.
We now aim at representing the discrete solution manifold by means of a second approximation process: the ROM approach.
The main goal is to rely on a low-dimensional space $Y_{N} \times Y_{N} \subset Y\disc \times Y\disc \subset Y \times Y$, with $N \ll N_h$, to solve a smaller system but preserving a certain accuracy in terms of errors between the reduced and the high-fidelity solution. The space $Y_{N}$ is a linear combination of state and adjoint \emph{snapshots}, i.e. high-fidelity solutions evaluated for properly chosen values of $\bmu \in \Cal P$. There are several strategies to build the reduced spaces: we postpone the description of them in the next sections and we now assume to be provided with these reduced spaces. Indeed, once built $Y_N \times Y_N$, we perform a Galerkin projection in the built low-dimensional framework to find the optimal solution for a new parametric instance: given $\bmu \in \mathcal P$, find $(y_N, p_N) \in Y_{N} \times Y_{N}$ such that
\begin{equation}
    \label{eq:weak_no_u_ROM}
    \begin{cases}
        m(y_N - y_{\text{d}}, w; \bmu) + a(w, p_N; \bmu) = 0                                              & \forall w \in Y_N, \\
        \displaystyle a(y_N, q; \bmu) - \frac{1}{\alpha}c(p_N, q; \bmu) =  \la G(\bmu), q \ra_{Y\dual, Y} & \forall q \in Y_N.
    \end{cases}
\end{equation}
The projection stage is convenient when it is independent from the high-fidelity dimension $N_h$.
To rely on a fast ROM solution, the approach should verify an \emph{offline-online decomposition} that consists in two stages:
\begin{itemize}
    \item[\small{$\circ$}] the \emph{offline phase}: the building process and the storing process. It depends on the dimension $N_h$ but it is performed only once;
    \item[\small{$\circ$}] the \emph{online phase}: given a new parameter, a Galerkin projection into the low-dimensional framework is performed to provide a solution in a small amount of time.
\end{itemize}
The aforementioned division is possible only when the problem verifies the \emph{affine decomposition}. In other words, when the weak forms  can be written as
\begin{equation}
    \label{eq:affine}
    \begin{matrix}
         & m(\cdot, \cdot; \bmu) =\displaystyle  \sum_{l=1}^{Q_m} \Theta_m^l(\boldsymbol{\mu})m^l(\cdot, \cdot), &
        \qquad  a(\cdot, \cdot; \bmu) =\displaystyle  \sum_{l=1}^{Q_a} \Theta_a^l(\boldsymbol{\mu})a^l(\cdot, \cdot),                                                                                                                                  \\
         & c(\cdot, \cdot; \bmu) =\displaystyle  \sum_{l=1}^{Q_c} \Theta_c^l(\boldsymbol{\mu})c^l(\cdot, \cdot), & \qquad  \la G(\boldsymbol{\mu}), q \ra =\displaystyle  \sum_{l=1}^{Q_{G}} \Theta_{G}^l(\boldsymbol{\mu})\la G^l, q \ra_{Y\dual, Y}, \\
    \end{matrix}
\end{equation}
for some $Q_m, Q_a, Q_c$ and $Q_G$ in $\mathbb N$, with $\Theta_m^l,\Theta_a^l,\Theta_c^l,$ and $\Theta_{G}^l$ smooth real functions depending on $\boldsymbol{\mu}$ and $m^l(\cdot, \cdot), a^l(\cdot, \cdot),c^l(\cdot, \cdot)$ and $ \la G^l, \cdot \ra_{Y\dual, Y}$ independent of $\bmu$.
The structure \eqref{eq:affine} is necessary to rely on an offline-online concept. 
In this way, in the offline stage we can assemble and store all the $\bmu$-independent quantities together with the reduced spaces basis functions, while, in the online phase, we compute the $\bmu$-dependent quantities for a specific parametric instance and the reduced optimality system is assembled and solved. \\
In our context, this assumption does not hold. 
Indeed, the bilinear form $c(\cdot, \cdot; \bmu)$ does not admit an affine decomposition, since the boundary integration over $\Gamma_C^{\mu_u}$ changes at each $\bmu$ and it has to be assembled for each parameter. \\
In order to recover the affine assumption, one can rely on several techniques, such as Empirical Interpolation Method (EIM) or Discrete Empirical Interpolation Method (DEIM). We refer the interested reader to \cite{barrault2004empirical, DEIM} and \cite[Chapter 5]{hesthaven2015certified}. We stress that in some of the results we are going to present in Section~\ref{sec:results}, we exploited DEIM to guarantee offline-online decomposition and, consequently, computational efficiency. We briefly recall DEIM process in Section~\ref{sec:DEIM}.

\subsection{The rational behind tailored reduced strategy for vbOCP($\boldsymbol \mu$)}
\label{sec:transport}

We observe from the numerical simulations of Section~\ref{sec:result_simple}, that the solution drastically changes with respect to $\mu_u$ and the problem might feature very complicated structures as transport-wave-like phenomena, i.e. moving fronts. 
The detail are described in Section~\ref{sec:results}, together with the test cases.

Namely, when dealing with vbOCP($\bmu$), standard reduced approaches may lead to inefficient reduced solution. Indeed, vbOCP($\bmu$)s need many basis functions to provide and accurate representation of the high-fidelity approximation. 
A large number of modes, combined with the non-affine structure of the problem, translates into unbearable reduced simulations, i.e.\ the reduced model is not competitive with respect to the high-fidelity one.

For this reason we decide to exploit strategies that are usually related to the field of nonlinear manifold reduction. 
In the last years, many approaches have been conceived to solve this kind of issues and much effort has been done to develop techniques capable to decrease the number of basis functions needed to efficiently reduce the system. They rely on pre-processed snapshots before applying the reduction algorithm  \cite{Iollo2014,Ohlberger2013901,nonino2019overcoming,PAPAPICCO2022114687,Rim2018118,Taddei2020A997} or on adapting the bases to the moving features \cite{Peherstorfer2020A2803, Zimmermann2018234}. Other approaches employ a transformation map of the solution in a reference domain to avoid the nonlinear features of the solution manifold, we refer the interested reader to \cite{Cagniart2019131, Ohlberger2013901, torlo2020model} and the reference therein. These strategies have been successfully employed for advection-dominated phenomena, most of all in a time-dependent setting.

Other strategies build local reduced basis since global structure may lead to inaccurate results if few modes are exploited. We postpone the discussion on local ROMs to Section~\ref{sec:LPOD}. Here, we report a far from exhaustive list of useful references on the topic \cite{Amsallem2016,amsallem2012nonlinear,Borggaard20162155,Dihlmann2012156,drohmann2011adaptive, Eftang2011395,Eftang2011179,Haasdonk2011423, Maday2013A2417}.

Moreover, we want to stress that a novel and successful research field focuses on the use of ROMs enhanced by artificial intelligence techniques to deal with nonlinear model order reduction, the interested reader may refer to \cite{Fresca2021, Fresca2022, Lee2020, romor2022non} and the reference therein.

\subsection{DEIM}
\label{sec:DEIM}
 The DEIM approach overcomes the problem of dealing with a non-affine system, as in our case. We follow the seminal paper \cite{DEIM}. We approximate the bilinear form $c(\cdot, \cdot, \boldsymbol{\mu})$ applying a projection onto a subspace of dimension
$N_{DEIM} \ll N_h$. We recall that the bilinear form related to the control is: 
$$
{c(u ,q; \bmu)} = \int_{\Gamma_{C}^{\mu_u}}u q \; \text{ds} = \int_{\Gamma_{C}^0} \chi_{\mu_u}uq\; \text{ds},
$$
where $\chi_{\mu_u} = 1$ where the control is applied and $\chi_{\mu_u} = 0$ elsewhere, while ${\Gamma_{C}^0} \subset \partial \Omega$ is the maximum portion of the boundary that can feature the control action. All the quantities are meant in the sense of traces. To apply the DEIM algorithm, we consider a discrete version of $\chi_{\mu_u}$, i.e. $\pmb{\chi}_{\mu_u} \in \mathbb R^{N_h}$: the FE coefficient array.
The goal is to approximate the characteristic function $\pmb{\chi}_{\mu_u}$ as
\begin{equation}
\label{eq:DEIMgoal}
\overline{\pmb{\chi}}_{\mu_u} = \sum_{q=1}^{N_{DEIM}} \Theta^q(\boldsymbol{\mu}) \mathsf z_q,
\end{equation}
with $\mathsf z_q \in \mathbb R^{N_h}$ basis vectors that do not depend on the parameter and $\Theta^q(\boldsymbol{\mu})$ are real coefficients for each value of the parameters. The basis vectors $\{\mathsf z_q\}_{q=1}^{N_{DEIM}}$ are built through a POD approach over the evaluation of several instances of $\pmb \chi_{\mu_u}$.

We call the basis matrix related to the DEIM approximation as 
$$
\mathsf Z^{DEIM} = [\mathsf z_1, \dots, \mathsf z_{N_{DEIM}}].
$$
Algebraically, \eqref{eq:DEIMgoal} can be written as 
$$\overline{\pmb \chi}_{\mu_u}= \mathsf Z^{DEIM}\mathsf \Theta(\boldsymbol{\mu}),$$ 
where $\mathsf \Theta(\boldsymbol{\mu}) = [\Theta^1(\boldsymbol{\mu}), \dots, \Theta^{N_{DEIM}}(\boldsymbol{\mu})]^T$ is the column vector of the coefficients. 
To complete the offline procedure we need to find a set of interpolation indices $\mathbb I \subseteq \{1, \dots, N_h\}$ which will be used in the online phase to find the specific coefficients $\Theta^q(\boldsymbol{\mu}^*)$ for a given $\boldsymbol{\mu}^* \in \mathcal P$. 
The set $\mathbb I$ is computed by the magic point algorithm, see e.g.\ \cite{EIM, magic}. 

In the online phase, for a specific parameter $\boldsymbol{\mu}^*$, we select the $N_{DEIM}$ rows given by the indices in $\mathbb I$ to determine the coefficient $\mathsf \Theta(\boldsymbol{\mu}^*)$. 
Thus, defining the matrix 
$$
\mathsf P = [\mathsf e^{\mathbb I_1}, \dots, \mathsf e^{\mathbb I_{N_{DEIM}}} ] \in \mathbb R^{N_h} \times \mathbb R^{N_{DEIM}},
$$
where $\mathsf e^{\mathbb I_i} \in \mathbb R^{N_h}$, $\mathsf e^{\mathbb I_i}_j = \delta_{ij}$ and $\delta_{ij}$ stays for the  Kronecker delta, the vector $\mathsf \Theta (\boldsymbol {\mu})$ can be found solving the following system:
$$
\mathsf P^T \pmb \chi_{\mu_u}= (\mathsf P^T \mathsf Z^{DEIM})\mathsf \Theta(\boldsymbol {\mu}).
$$
Namely, 
$$\overline{\pmb \chi}_{\mu_u}  = \mathsf Z^{DEIM}\mathsf \Theta(\boldsymbol{\mu}) = \mathsf Z^{DEIM}(\mathsf P^T \mathsf Z^{DEIM})^{-1} \mathsf P^T \pmb \chi_{\mu_u} .$$
Once we are provided of this representation, the affine structure of the problem is recovered and the online-offline paradigm holds.\\
In the numerical results of Section~\ref{sec:results}, we will see the benefits of enhancing ROMs by DEIM technique.
\subsection{POD}
\label{sec:POD}
We focus on the building process based on POD-Galerkin strategy. The interested reader may refer to \cite{burkardt2006pod, Chapelle2013, hesthaven2015certified} for a detailed discussion on this algorithm. We here summarize the main features of it together with its generalization to OCP($\bmu$)s (that it is suited for vbOCP($\boldsymbol \mu$), too). 
We start sampling $(\bmu_1, \dots, \bmu_{N_{\text{max}}}) \in \mathcal P$ and computing the related high-fidelity approximations $$(y^{N_h}(\bmu_1), \dots, y^{N_h}( \bmu_{N_{\text{max}}}) \text{ and } (p^{N_h}(\bmu_1), \dots, p^{N_h}( \bmu_{N_{\text{max}}}).$$
From now on, we describe the process for the state variable only, since we apply a separate POD for the two variables: the procedure repeats for the adjoint variable, analogously.
\\ We aim at building $N < N_{\text{max}}$ basis functions by means of snapshots manipulation, to discard redundant information from the sampled parametric solutions. \\
We assume that the number of samples $N_{\text{max}}$ is large enough to well represent the solution manifold \eqref{eq:solution_manifold}. The process provides the $N-$dimensional reduced space $Y_N^y$ (and not $Y_N$) that minimizes the following quantities:
\begin{equation*}
    \sqrt{\frac{1}{N_{\text{max}}}
    \sum_{i=1}^{N_{\text{max}}} \min_{\zeta_N \in {Y}_{N}^y}\norm{y\disc(\boldsymbol{\mu}_i) - \zeta_N}_{Y}^2.
    }
\end{equation*}
To reach the minimization goal, we define the correlation matrix $\mathbf C^{y } \in \mathbb R^{N_{\text{max}} \times N_{\text{max}}}$ of state snapshots and we solve the eigenvalue problem
$
    \mbf C^{y} \omega_n^{y } = \lambda_n^{y } \omega_n^{y }
$ for $ 1 \leq n \leq N_{\text{max}},$ with $\norm {\omega_n^{y }}_{{Y}} = 1$. Due to the definition of correlation matrix, we can order the all-positive eigenvalues as $\lambda_1^{y } >\dots > \lambda_{N_{\text{max}}}^{y }> 0$ and retain the first $N$ eigenpairs $(\lambda_n^y, \omega_n^y)$ for $1 \leq n \leq N$. Finally, the basis are built as \cite{hesthaven2015certified, quarteroni2015reduced}, i.e.\
\begin{equation}
    \label{eq:basis}
    \xi_n^{y } = \displaystyle \frac{1}{\sqrt{{\lambda_n^{y }}}}\sum_{m = 1}^{{N_{\text{max}}}} (\omega_n^{y })_m y\disc (\boldsymbol{\mu}_m), \hspace{1cm} 1 \leq n \leq N.
\end{equation}
The choice of $N_{\text{max}}$ and $N$ can be made by studying the behaviour of $\lambda_n^y$ for $1 \leq n \leq N_{\text{max}}$. Indeed, defining as  $P_N: Y \rightarrow Y_N^y$ the projector from $Y$ onto $ Y_N^y$, the following relation holds \cite{hesthaven2015certified, quarteroni2015reduced}:
\begin{equation}
    \sqrt{\frac{1}{N_{\text{max}}}
    \sum_{i = 1}^{N_{\text{max}}}  \norm{y\disc(\boldsymbol{\mu}_{i}) - P_N(y\disc(\boldsymbol{\mu}_i)) }_{Y}^2} = \sqrt{
    \sum_{i = N + 1}^{N_{\text{max}}}\lambda_m^y.}
\end{equation}
Namely, a fast decay of the eigenvalue magnitude guaratees a good representation of the high-fidelity solution with a few basis functions.\\
The application of the POD to the adjoint variable leads to another reduced space $Y_N^p$. In principle, it can be different from $Y_N^y$ and this does not guarantee the well-posedness of the problem as proposed in Section~\ref{sec:hf}. Thus, we need a post-processing step over the basis. We apply the \emph{aggregated space} technique to build a common space for state and adjoint, a standard approach when dealing with OCP($\bmu$)s, see e.g.\ \cite{bader2016certified,bader2015certified,dede2010reduced,gerner2012certified,karcher2014certified,karcher2018certified, negri2015reduced,negri2013reduced,quarteroni2007reduced} as references. \\
The final reduced space is of the form
$$
    Y_N = \{(\xi^{y}_n, \xi^{p}_n), \; n = 1, \dots, N\}.
$$
In this way, we can define the basis matrix
$\mathsf Q = [\xi_{1}^{y} , \cdots , \xi_{N}^{y}, \xi_{1}^{p} ,\cdots , \xi_{N}^{p}] \in \mathbb R^{N_{h} \times 2N}$ and employ a Galerkin projection over the high-fidelity quantities of \eqref{eq:one_shot} we stored in the offline phase, dealing with a ${4N \times 4N}$ system of the form
\begin{equation}
    \label{eq:projN}
    \begin{bmatrix}
        \mathsf M_N & \mathsf D_N^T \\
        \mathsf D_N & \mathsf C_N           \\
    \end{bmatrix}
    \begin{bmatrix}
        \mathsf y_N \\
        \mathsf p_N
    \end{bmatrix} =
    \begin{bmatrix}
        \mathsf {y_{\text{d}}}_N \\
        \mathsf f_N              \\
    \end{bmatrix},
\end{equation}
where
$\mathsf M_N = \mathsf Q ^T M_{\text{o}}(\bmu) \mathsf Q$, $\mathsf D_N = \mathsf Q^T  \mathsf D_a(\bmu) \mathsf Q$,   $\mathsf C_N = \mathsf Q^T  \mathsf C(\bmu) \mathsf Q$, $\mathsf {y_{\text{d}}}_N  = \mathsf Q^T  \mathsf {y_{\text{d}}}_N$ and $\mathsf f_N = \mathsf Q^T  \mathsf f$.
We stress that if the affine decomposition holds, the high-fidelity matrices inherit the affine property from their continuous counterpart \eqref{eq:affine} and the reduced $\bmu$-independent part of the matrices can be stored once and for all in the offline phase.

\subsection{L-POD}
\label{sec:LPOD}
To increase the accuracy of the reduced order model, the sampled solution manifold $\mathbb M^{N_h}$ is divided into different \emph{snapshot regions} and for each region a reduced model is built. The idea of local POD is not new in literature. Indeed, building separate sets of basis functions for a subspace of the parametric solution manifold turns out to be beneficial to reduce the dimensionality of complicated problems. The interested reader may refer to several papers, such as \cite{Amsallem2016,amsallem2012nonlinear,Borggaard20162155,Dihlmann2012156,drohmann2011adaptive, Eftang2011395,Eftang2011179,Haasdonk2011423, Maday2013A2417}. Local reduced order modelling has been successfully employed in many fields of applications, for example in cardiac models \cite{Pagani, VenezianiLocal}, computational fluid dynamics and aerodynamics \cite{washabaugh2012nonlinear}.
The accuracy of the local bases is higher with respect to a global approach when one is dealing with nonlinear and non-affine problems.
Let us assume to have clustered the $N_{max}$ snapshots (we postpone the discussion on \emph{how cluster them} later on) in $J$ groups. In this way, we divide $\mathbb M\disc$ in sub-regions $\{\mathbb M^{N_h}_i\}_{i=1}^J$ and we want to separately reconstruct the local solution manifolds $\mathbb M^{N_h}_i$ for $i = 1, \dots, J$. The new offline phase consists in the application of the POD strategy enhanced with aggregated spaces for each of the sub-manifold. This translates in $J$ basis matrices $\mathsf Q_i$ for $i = 1, \dots, J$, that are able to locally represent the system. Indeed, the online procedure first sorts a specific parameter, say $\boldsymbol {\mu}^*$, into a group and, then, applies the specific POD-Galerkin procedure over the selected local space.

\begin{algorithm}[H]
\label{alg1}
\caption{Pseudo-code for L-POD (Offline phase)}\label{al:01}
\begin{algorithmic}[1]
\State{$N_{\text{max}}, N, \tau,$ LIST $= [I_{\mu_u}], M, i=0,  [\bmu_{1}, \cdots, \bmu_{N_{\text{max}}}], J=0$}\Comment{Inputs}
\For{$\bmu \in [\bmu_{1}, \cdots, \bmu_{N_{\text{max}}}]$} 
\State{Solve \eqref{eq:one_shot} and store $y(\bmu)$ and $p(\bmu)$}
\EndFor
\While{ LIST $\neq \emptyset$ or $i < M$}
\For {$I \in LIST$}
\State{LIST = LIST$\setminus I$}
\For {$\mu_u \in  [{\mu_u}_{1}, \cdots, {\mu_u}_{N_{\text{max}}}]$}
\If {$\mu_u \in I$}
\State{Enrich local state and adjoint covariance matrices}
\EndIf
\EndFor
\State{Solve the $N-$eigenvalue problem over the local covariance matrices}
\If{$\lambda_N^y > \tau$}
 \State{$i = i +1$}
\State{Halve $I$ in $I_1$ and $I_2$}
\State{LIST = LIST $ \cup \; I_1 \cup I_2$}
\Else
\State{$J = J+1$}
\State{Build the local $N-$dimensional basis (aggragated spaces)}
\State{Build the local reduced operators}
\EndIf
\EndFor
\EndWhile
\end{algorithmic}
\end{algorithm}
It remains to understand how to cluster the snapshots and how to sort the online parameters. 
In literature, many strategies have been employed. For time dependent problems, a classical approach is to divide the considered time interval in time windows and build a reduced space for each one of them, as done in \cite{Borggaard20162155, Dihlmann2012156,drohmann2011adaptive}. Another strategy relies on the division of the parametric space $\mathcal P$ adaptively, see e.g.\ \cite{Eftang2011395,Eftang2011179,Haasdonk2011423}. Last, one can use classification algorithms to cluster the snapshots in a beneficial way, as done in \cite{Pagani}.\\
We propose a tailored algorithm for this specific problems, guided by the numerical results we are going to show in Section~\ref{sec:results}.\\
As we already specified in Section~\ref{sec:transport}, the change of $\mu_u$ plays a crucial role in the complexity of the problem. Thus, our choice was to adaptively divide the geometrical parameter interval $ I_{\mu_u}$, where $\mu_u$ lives. First of all, we fix a tolerance $\tau$, a maximum value of divisions $M$ and a basis number $N$. We applied a first POD on the whole interval $I_{\mu_u}$ for a given basis number $N$. If the state eigenvalue $\lambda_N^y \geq \tau$, we halve $I_{\mu_u}$ in  $I_{\mu_u}^1$ and $I_{\mu_u}^2$. We proceed analogously for the two sub-intervals. If the criterion is verified we stop. If not, we continue halving once again the interval we are dealing with. The procedure ends when all the intervals meet the criterion, or when $I_{\mu_u}$ is divided in $M$ sub-intervals. This procedure leads to a partition of the geometrical parametric interval as
\begin{equation*}
I_{\mu_u} = \bigcup_{i=1}^J I_{\mu_u}^i,
\end{equation*}
and to the creation of $J$ reduced spaces, one for each of the sub-intervals, with $J \leq M$.
{Summing up, in the offline phase:
\begin{itemize}
\item[$\circ$] we generate all the high-fidelity solutions and we apply a first POD;
\item[$\circ$] if the criterion over $\lambda_N^y$ is not verified, we halve $I_{\mu_u}$ and we separate the snapshots, consequently. We perform the POD operation on each sub-interval. We iteratively repeat this process;
\item[$\circ$] when we meet the criterion over $\lambda_N^y$ for a sub-interval, we build and store the local basis functions. 
\end{itemize}
If the tolerance on the state eigenvalues is not reached, the procedure stops after a maximum number of iterations.}\\
In the online phase,  it suffices to sort $\boldsymbol {\bmu}$ with respect to the sub-interval it belongs to.
We recap the offline procedure in Algorithm \ref{al:01}. We took inspiration from the time-windowed strategies presented in \cite{Borggaard20162155, Dihlmann2012156,drohmann2011adaptive} adapting them to the vbOCP($\boldsymbol \mu$). We recall that the partition of $I_{\mu_u}$ looks a natural choice, since the transport phenomena are directly related to $\mu_u$, as we will discuss in Section~\ref{sec:results}.
{
\begin{remark}
The strategy we propose is only one path one can take to build local spaces. We believe that many different adaptive strategies, even more efficient, can be found to deal with vbOCP($\boldsymbol \mu$), as we stressed out in the beginning of this section.
However, our goal is to underline the main features of the vbOCP($\boldsymbol \mu$) model as a first step towards a more complete and deeper analysis. Thus, we did not try other techniques since such an experimental study is beyond the scope of the contribution.\\
Furthermore, we stress that the adaptive strategy is only related to the geometrical parameter interval $I_{\mu_u}$. The extension to three-dimensional (3D) problem is natural, since the vbOCP($\boldsymbol \mu$) can be represented by a characteristic function $\chi_{\mu_u}(x_1, x_2, x_3)$, where $x_i$ for $i=1,2,3$ denotes the coordinates of the domain $\Omega$.
\end{remark}
}

\subsection{Geo-R}
\label{sec:geo-r}
When dealing with simple settings, it is possible to recast the problem as formulated in \cite{RozzaHuynhPatera2008}. 
The strategy, in the standard OCP($\bmu$) framework, was already exploited in several works, see e.g.\ \cite{StrazzulloZuazua,negri2015reduced,negri2013reduced,Strazzullo2,StrazzulloRB}.
The aim is to solve the optimality FE systems \eqref{eq:weak_no_u_FE} and the reduced one \eqref{eq:weak_no_u_ROM} in a reference domain $\Omega_o \eqdot \Omega(\mu_u^o)$. In our case, this means to recast the problem in a framework where the control boundary and, consequently, the Neumann boundary, are ``fixed", simplifying the setting at hand.
{In the following, we indicate with $\mathring{\Xi}$ and $\overline{\Xi}$ the internal part and the closure of a spatial domain $\Xi$, respectively. } \\
We assume that the domain $\Omega_o$ is the union of $R$ non-overlapping subdomains $\Omega^r_o$, i.e.\
\begin{equation}
    \Omega_o = \bigcup_{r=1}^R \Omega^r_o \text{ with } \mathring{\Omega}^{r'}_o \cap \mathring{\Omega}^{\bar r}_o = \emptyset \text{ for } 1 \leq r' < \bar r \leq R.
\end{equation}
The same assumption holds true for the domain $\Omega$, i.e.
\begin{equation}
    \Omega = \bigcup_{r=1}^R \Omega^r\text{ with } \mathring{\Omega}^{r'} \cap \mathring{\Omega}^{\bar r} = \emptyset \text{ for } 1 \leq r' < \bar r \leq R.
\end{equation}
We want to define {a piece-wise} affine map $T_{\mu_u}: \Omega_o \rightarrow \Omega$ such that $\Omega = T_{\mu_u}(\Omega_o)$. The map $T_{\mu_u}$ is {made by local maps $T_{\mu_u}^r$ defined over the subdomains $\Omega^r_o$ and on each subdomain $T_{\mu_u}^r$ is invertible. We remark the $T_{\mu_u}^r(\Omega_o^r) = \Omega^r$}. Moreover, the local maps glue continuously, namely, given a point $\mathbf x \in \Omega_o$:
{
\begin{equation}
T_{\mu_u}^{r'}(\mathbf x_o) = T_{\mu_u}^{\bar r}(\mathbf x_o) \quad \forall \mathbf x_o \in
    \overline{\Omega}^{r'}_o \cap \overline{\Omega}^{\bar r}_o, \text{ for } 1 \leq r' < \bar r \leq R.
\end{equation}}
{The last assumption is that the local maps are affine transformations}. Indeed, for every $r=1,\dots, R$ we require
${T_{\mu_u}^r(\mathbf x_o)} = \mathsf c^r(\mu_u) + \mathsf G^r(\mu_u) \mathbf x_o$, where $\mathsf c^r(\mu_u)\in \mathbb R^2$ is a translation vector and $\mathsf G^r(\mu_u) \in \mathbb R^{2\times 2}$ is a linear transformation matrix. Once defined the map, it is clear that the continuous version of the optimality system \eqref{eq:weak_no_u}, together with the FE one \eqref{eq:weak_no_u_FE} and the reduced one \eqref{eq:weak_no_u_ROM}, can be solved in the new reference domain after a change of variables over the integral forms using the inverse { of the local maps}. 
We refer the reader to \cite{RozzaHuynhPatera2008} for the details. As already specified, the main advantage is that {we work with a reference $\Gamma_C^{\mu_u^o}$}, i.e.\ all the snapshots are taken in a new framework where the control boundary stays unchanged. Once the change of variable is performed, one can proceed with standard ROM techniques, as POD.
Once the map is found, after the change of variable problem \eqref{eq:weak_no_u} becomes: given $\bmu \in \mathcal P$, find $(y,p) \in Y\times Y$ such that
\begin{equation}
\label{eq:weak_geo}
\begin{cases}
\overline{m}(y - y_{\text{d}}, w; \bmu) + \overline{a}(w, p; \bmu) = 0 & \forall w \in Y,\\
\displaystyle \overline{a}(y, q; \bmu) - \frac{1}{\alpha}\overline{c}(p, q; \bmu) =  \la \overline{G}(\bmu), q \ra_{Y\dual, Y}& \forall q \in Y,
\end{cases}
\end{equation}
where $Y=H^1_0(\Omega_o)$ is this framework.

\section{Numerical Results}
\label{sec:results}
In this section we present two different tests over different geometries. The governing PDE($\bmu$) is an advection-diffusion equation for both the experiments. The problem reads:
given $\bmu \in \Cal P$, find the pair $(y,p) \in Y \times Y$ such that
\begin{equation}
  \label{eq:test_1}
  \begin{cases}
    \displaystyle y \chi_{\Omega_{\text{obs}}} - \frac{1}{\mu_1}\Delta p - x_2(1 -x_2)\frac{\partial p}{\partial x_1} = 
    \mu_2 \chi_{\Omega_{\text{obs}}}                                                           & \text{ in } \Omega,              \\
    \displaystyle  - \frac{1}{\mu_1}\Delta y + x_2(1 -x_2)\frac{\partial y}{\partial x_1} =  0 & \text{ in } \Omega, \vspace{2mm} \\
    \displaystyle \frac{1}{\mu_1}\frac{\partial y}{\partial n} = 0                             & \text{ on $\Gamma_N^{\mu_u},$}   \\
 \displaystyle \frac{1}{\mu_1}\frac{\partial p}{\partial n} +  x_2(1 -x_2) n_1 p= 0                             & \text{ on $\Gamma_N^{\mu_u},$}   \\
    \displaystyle \frac{1}{\mu_1}\frac{\partial y}{\partial n} =  \frac{1}{\alpha}p            & \text{ on $\Gamma_C^{\mu_u},$}   \\
    y = 1 \text{ and } p = 0                                                                   & \text{ on $\Gamma_{D},$}         \\
  \end{cases}
\end{equation}
where $(x_1, x_2)$ denotes the spatial coordinates of the domain $\Omega$ and $n_1$ is the first component of the outer vector to $\Gamma_N^{\mu_u}$. We also refer to the vector of spatial coordinates as $\mathbf x \eqdot [x_1, x_2]^T$.\\ We recall that, despite the same constraint, the geometry of the test cases drastically changes. This is the key aspect of the numerical results we are going to show: for standard and simple geometries the problem simplifies and reference-domain-based techniques can be exploited in order to achieve satisfactory results. This is not the case for complex geometries (nonlinear boundaries, holes...), where these techniques can be hardly applied.
\subsection{Test Case 1: POD, L-POD and Geo-R over a standard geometry}
\label{sec:result_simple}
\begin{figure}
  \centering
  \includegraphics[width=\textwidth]{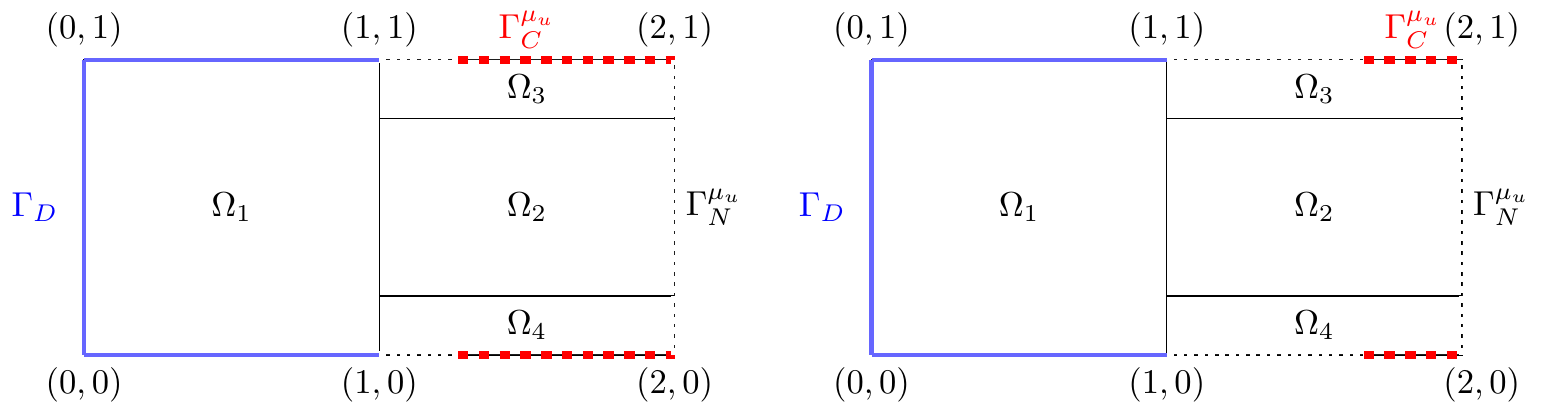}
  \caption{(Test 1). Domain $\Omega$. \textit{Observation domain:} $\Omega_{\text{obs}} = \Omega_3\cup \Omega_4$, \textit{Control domain:} $\Gamma_C^{\mu_u}$ (red dashed line). \textit{Blue solid line:} Dirichlet boundary conditions. \textit{Black dotted line:} Neumann boundary conditions. We represent the domain for $\mu_u = 0.3$ (left) and $\mu_u = 0.7$ (right).}
  \label{fig:dominio_1}
\end{figure}

We propose a vbOCP($\bmu$) governed by \eqref{eq:test_1} with two physical parameters and the control boundary that changes with respect to $\mu_u$, as depicted in Figure~\ref{fig:dominio_1}.  The observation domain is $\Omega_{\text{obs}} = \Omega_3 \cup \Omega_4$ where $\Omega_3 = [1, 2]\times [0.8, 1]$,  $\Omega_4 = [1, 2]\times [0, 0.2]$, while $\Omega_1$ is the unit square and $\Omega_2 = [1, 2]\times [0.2, 0.8]$. The control acts on the boundary $\Gamma_C^{\mu_u} =([1+\mu_u, 2] \times \{0\}) \cup  ([1+\mu_u, 2] \times \{1\}).$ Also the portion of the domain where Neumann boundary conditions are applied changes with respect to $\mu_u$: $\Gamma_N^{\mu_u} = ([1, 1+\mu_u) \times \{0\}) \cup ([1, 1+\mu_u) \times \{1\})\cup (\{2\}\times [0,2])$.

The parameter is $\bmu \eqdot (\mu_1, \mu_2, \mu_u) \in
  \Cal P = (6.0, 20.0)  \times (0.5, 3.0)\times (0.0, 1.0)$, where $\mu_1$ is related to the P\'eclet number, $\mu_2$ describes the constant desired state we want to reach in $\Omega_{\text{obs}}$ and $\mu_u$, as already specified, changes the boundary portion where the control acts. In analogy with the numerical results shown in \cite{karcher2018certified,negri2013reduced,Strazzullo2,StrazzulloRB}, we choose $\alpha = 0.07$.
We now describe the performances of standard POD.

\begin{figure}
  \centering
  \includegraphics[width=0.48\textwidth]{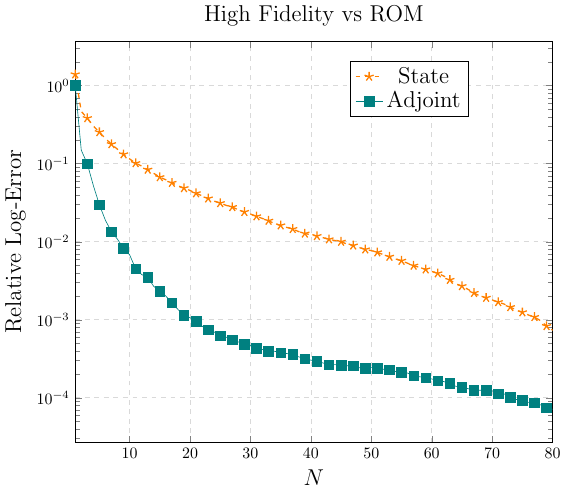}\
  \includegraphics[width=0.48\textwidth]{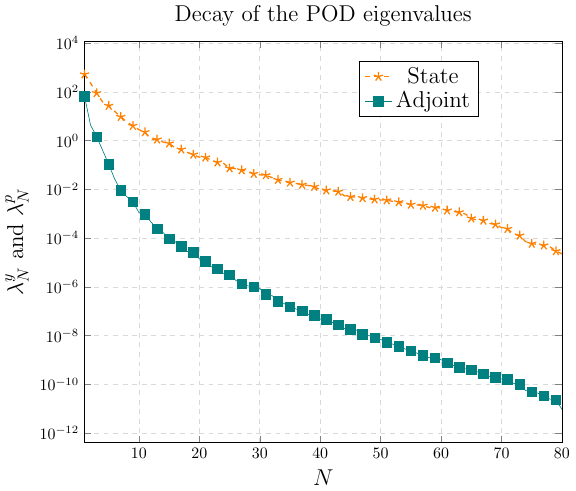}
  \caption{(Test 1: POD). {\emph{Left}. Averaged relative log-error for the two variables. \emph{Right}. Decay of the eigenvalues for the two variables}.}
  \label{fig:error_lambda300}
\end{figure}

First of all, we recover the affine decomposition through a DEIM approach applied to a training set of $350$ uniformly distributed parameters. The method stops when the DEIM approximation is a good representation of the actual problem (the stopping criterion over the eigenvalue decay is set to be lower then $10^{-5}$  and leads to $N_{DEIM}=80$). Thus, a standard POD is performed on other $300$ snapshots picked by means of uniform distribution. The high-fidelity dimension of the $\mathbb{P}^1-\mathbb{P}^1$ approximation is $2N_h = 10366$. In Figure~\ref{fig:error_lambda300} (left plot) we show the averaged relative errors for the state and the adjoint variable between the high-fidelity and reduced solutions. We plot the average value over a testing set of 150 uniformly distributed parameters in $\mathcal P$ of the quantities
$$
  e_y = \frac{\norm{y^{N_h} - y_N}_{Y}}{\norm{y^{N_h}}_Y} \quad \text{and}\quad e_p = \frac{\norm{p^{N_h} - p_N}_{Y}}{\norm{p^{N_h}}_Y}.
$$

\begin{figure}
  \centering
  \includegraphics[width=0.3\textwidth]{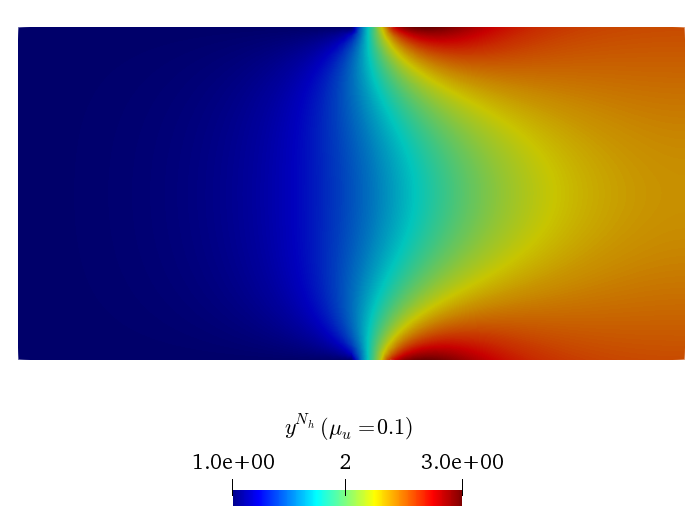}
  \includegraphics[width=0.3\textwidth]{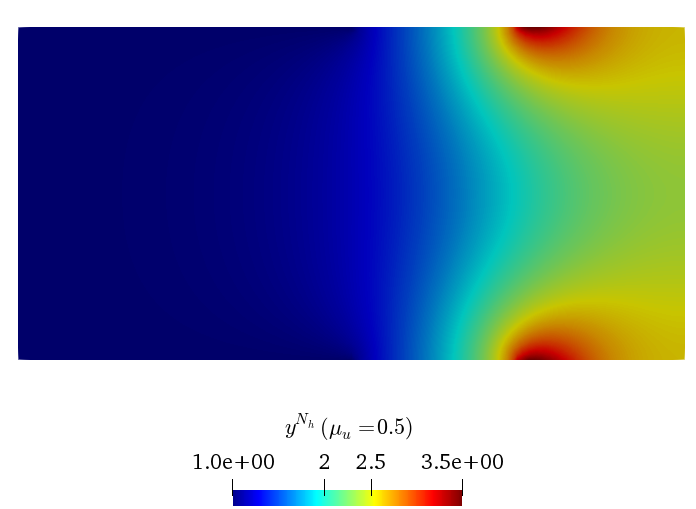}
  \includegraphics[width=0.3\textwidth]{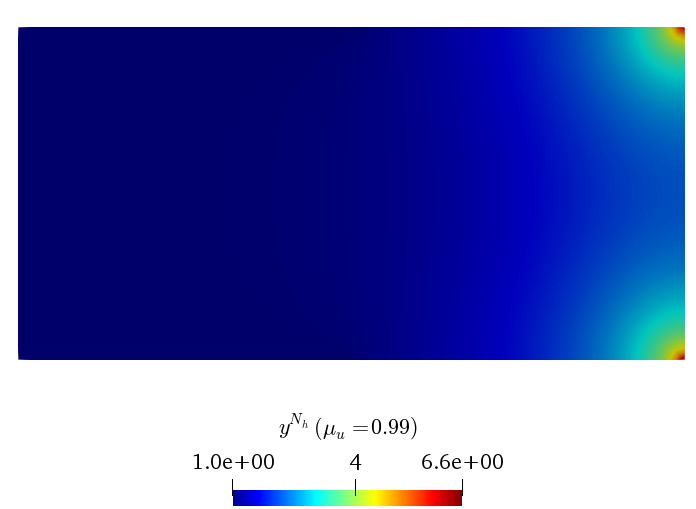}\\
  \includegraphics[width=0.3\textwidth]{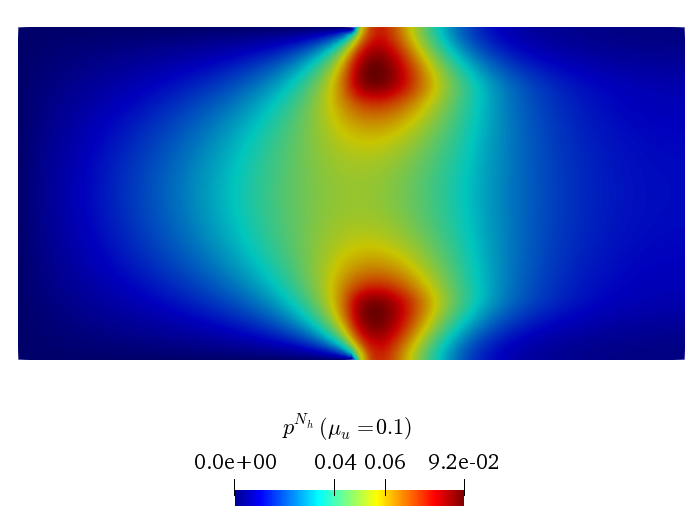}
  \includegraphics[width=0.3\textwidth]{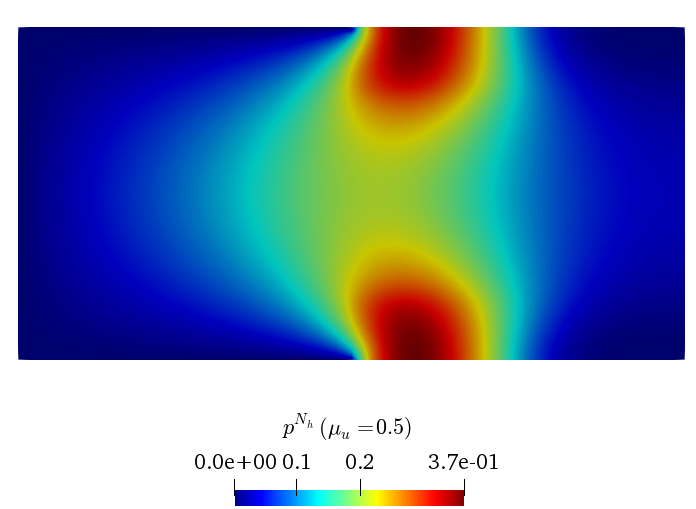}
  \includegraphics[width=0.3\textwidth]{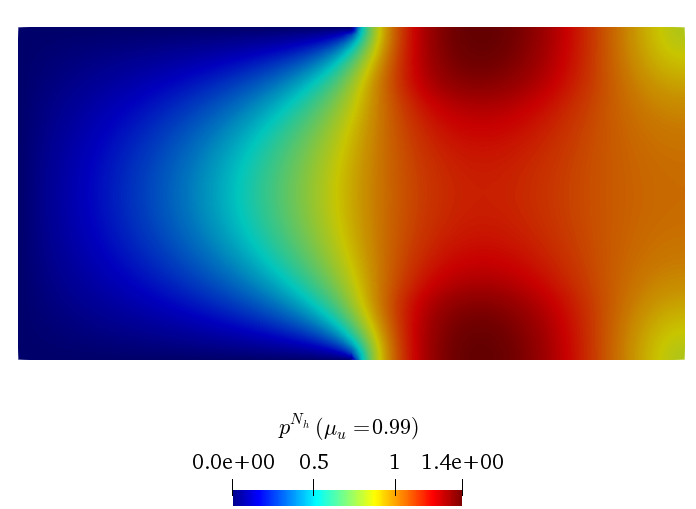}\\
  \caption{(Test 1: POD). Some high-fidelity solutions for fixed physical parameters ($\mu_1 = 12, \mu_2 = 2.5$) and varying $\mu_u = 0.1, 0.5, 0.99$, from left to right. The state solutions and the adjoint ones are represented on the top and on the bottom of the Figure, respectively.}
  \label{fig:snap}
\end{figure}

\begin{figure}
  \centering
  \includegraphics[width=0.48\textwidth]{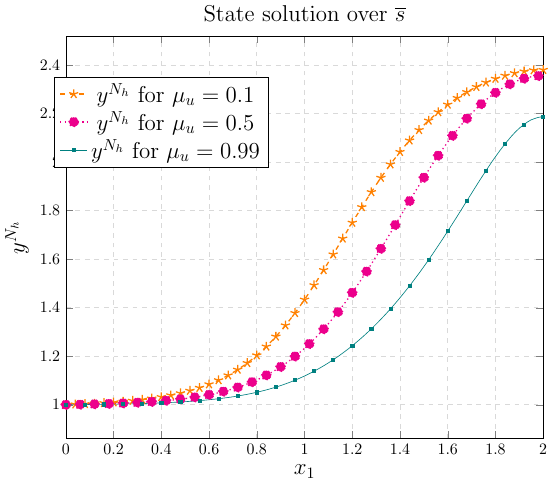}\
  \includegraphics[width=0.48\textwidth]{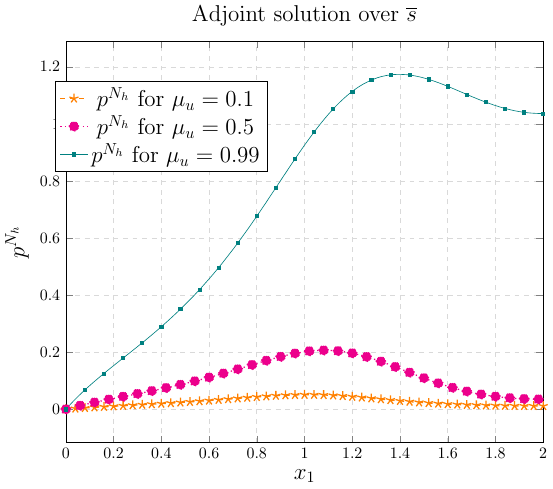}
  \caption{(Test 1: POD). {\emph{Left}. State solution $y^{N_h}$ for fixed physical parameters $\mu_1=12$ and $\mu_2=2.5$ and varying $\mu_u \in \{0.1, 0.5, 0.99\}$ over the segment $\overline s = [0,2] \times \{0.5\}$. \emph{Right}. Analogous representation for the adjoint solution}.}
  \label{fig:snap_over_line}
\end{figure}

We denote these averaged quantities with $E_y$ and $E_p$, respectively.
It is clear that, too many basis functions are needed in order to reach acceptable values for $E_y$, say $N = 80$ to reach $E_y \sim O(10^{-3})$.
We remark that $N = 80$ translates in a reduced system of dimension $4N \times 4N$ to be solved, by means of the aggregated space technique. The adjoint error $E_p$ takes only $N=20$ to reach such a threshold. This behaviour is not unexpected if one looks at the eigenvalue decay related to the POD modes, depicted in Figure~\ref{fig:error_lambda300} (right plot). It is clear that the POD modes are struggling in representing the state variable, that has a slower decay of the eigenvalues with respect to the adjoint variable. The reason is clearer if we look at Figure~\ref{fig:snap}. We fix the physical parameters to $\mu_1 = 12$ and $\mu_2 = 2.5$ and we study the behaviour of the solution with respect the ``vanishing'' action of the control, represented in this case by $\mu_u \in \{0.1, 0.5, 0.99\}$, from left to right. The top plots describe the state solution $y^{N_h}$. Here, we can observe that changing $\Gamma_C^{\mu_u}$ causes a \emph{transportation} of the temperature fields. On the upper and lower boundaries, we see a peak due to the action of the control that heats the system. As one can see, the peak magnitude increases and shifts from left to right. The transport issue is verified by the left plot of Figure~\ref{fig:snap_over_line}, too: it represents $y^{N_h}$ over the segment $\overline s = [0,2] \times \{0.5\}$. Also in this case we fix $\mu_1 = 12$ and $\mu_2 = 2.5$ and $\mu_u \in \{0.1, 0.5. 0.99\}$. From the plot we see the solution moving and transporting itself along the $x_1-$axis almost unchanged. This behaviour is expected for all the points in $\Omega_2 \cup \Omega_{\text{obs}}$. This is not the case for the adjoint solution (compare the bottom plots of Figure~\ref{fig:snap} and the right plot of Figure~\ref{fig:snap_over_line}).

Namely, letting the control boundary $\Gamma_C^{\mu_u}$ change makes the problem a very difficult task to tackle with standard POD (and in general with ROM, as we specified in Section~\ref{sec:transport}). This is confirmed also from the fact that in \cite{karcher2018certified,negri2013reduced,Strazzullo2,StrazzulloRB}, works where the control boundary was fixed, a few basis functions were necessary to represent the whole phenomenon. The real issue is the moving boundary.
We tried two ways to overcome this problem: (i) the L-POD and (ii) the Geo-R.
In the following we compare the results of the approaches.

\subsubsection{L-POD}
\label{subsub:LPOD}
The L-POD approach, as specified in Section~\ref{sec:LPOD}, is capable to reach more accurate results with respect to the high-fidelity approximation. The setting is the same proposed for the standard POD. We recall that $\bmu \eqdot (\mu_1, \mu_2, \mu_u) \in
  \Cal P = (6.0, 20.0)  \times (0.5, 3.0)\times (0.0, 1.0)$ and $\alpha = 0.07$.
Also in this case, the affine decomposition is recovered by means of a DEIM strategy applied over $350$ uniformly distributed parameters. The DEIM is applied with a tolerance of $10^{-5}$ over the eigenvalue decay and leads to $79$ basis functions (this number may not coincide with the standard POD since the $350$ parameters are different with respect to the ones used in the standard POD). We stress that the DEIM is applied globally, for all the parametric space $\mathcal P$. We exploit Algorithm \ref{al:01}. We set the maximum number of iteration $M = 10$ and $\tau = 10^{-3}$ with $N = 30$. We computed $300$ snapshots and they are collected in the $J = 4$ groups determined by the adaptive algorithm. 
Namely, we did not reach the value of $M$ and $I_{\mu_u}$ is partitioned as follows:
$$I_{\mu_u} = (0, 0.25] \cup (0.25, 0.5] \cup (0.5, 0.75] \cup (0.75,1).$$
{Thanks to the application of four separate PODs, we were able to provide (locally), a faster decay of the eigenvalues.}

\begin{figure}
  \centering
  \includegraphics[width=0.48\textwidth]{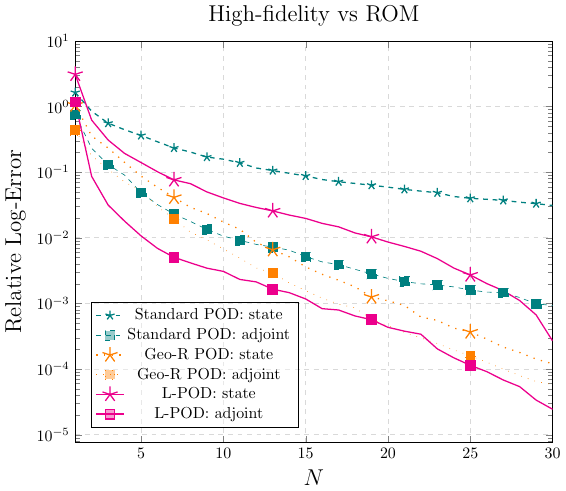}\
  \includegraphics[width=0.48\textwidth]{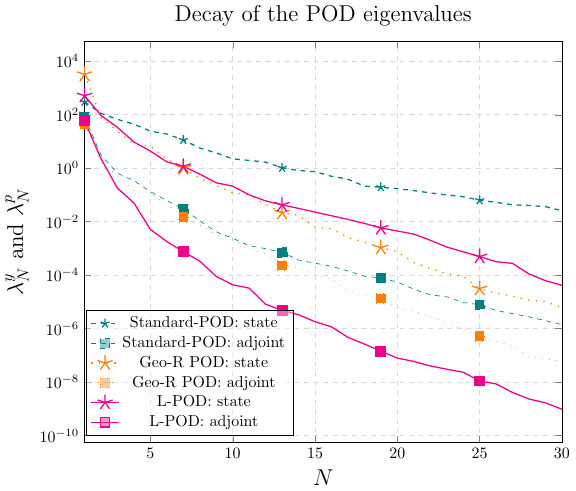}
  \caption{(Test 1: approaches comparison). {\emph{Left}. Averaged relative log-error for the two variables for all the approaches. \emph{Right}. Decay of the eigenvalues for the two variables for all the approaches}.}
  \label{fig:error_comp}
\end{figure}

{In the right plot of Figure~\ref{fig:error_comp}, we show the decay of the averaged eigenvalues $\bar \lambda_n^y$ and $\bar \lambda_n^p$ for $n=1,\dots, N$, defined, in the setting of L-POD, as
\begin{equation}
  \label{eq:eigLPOD}
  \bar \lambda_n^y = \frac{1}{J} \sum_{j = 1}^J {\lambda_n^y}^j \text{ and } \bar \lambda_n^p = \frac{1}{J} \sum_{j = 1}^J {\lambda_n^p}^j.
\end{equation}
For $N=30$, the L-POD eigenvalues are two orders of magnitude smaller than the ones related to standard POD, both for state and adjoint variables. Finally, the state eigenvalues decay is faster for Geo-R with respect to L-POD, while the adjoint eigenvalues have an opposite behaviour.
}
This reflects onto the average error over $150$ uniformly distributed parameters depicted in the left plot of Figure~\ref{fig:error_comp}. Each $\bmu$ in the online phase is sorted in the respective sub-interval and the respective POD basis are used for the Galerkin projection. This local reconstruction gives very accurate results: for $N=30$, $E_y \sim 2\cdot 10^{-3}$ and $E_p \sim 2\cdot 10^{-4}$. We stress that the standard POD is not capable to reach these values neither with $N = 80$. {We postpone the comparison between the L-POD and Geo-R relative errors in the next section.}

\subsubsection{Geo-R}
\label{subsub:geor}
As described in Section~\ref{sec:geo-r}, we are going to solve the vbOCP($\boldsymbol \mu$) in a reference domain $\Omega_o$.
We recall that in the examples usually proposed in literature, the control was ``fixed'' and the geometrical parameter did not change its nature as in our case.

\begin{figure}
  \centering
  \includegraphics[width=0.5\textwidth]{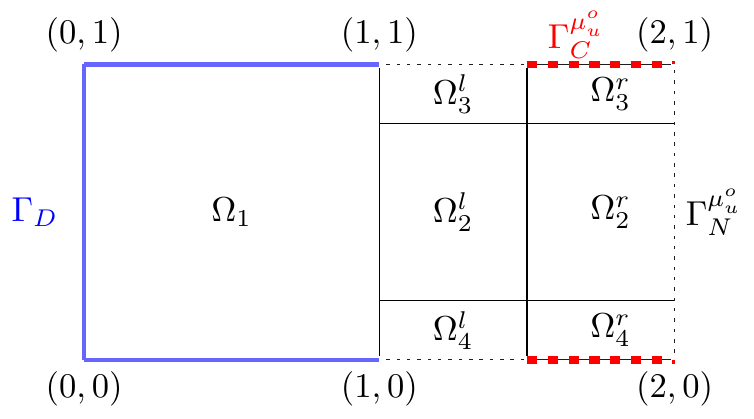}
  \caption{(Test 1). Domain $\Omega_o$ for $\mu_u^o = 0.5$. \textit{Observation domain:} $\Omega_{\text{obs}} = \Omega_3^l\cup  \Omega_3^r\cup \Omega_4^l\cup  \Omega_4^r$, \textit{Control domain:} $\Gamma_C^{\mu_u^o}$ (red dashed line). \textit{Blue solid line:} Dirichlet boundary conditions. \textit{Black dotted line:} Neumann boundary conditions.}
  \label{fig:dominio_recast}
\end{figure}

In Figure~\ref{fig:dominio_recast} we show the geometrical domain used in the Geo-R framework. 
Here, we consider, once again, $\bmu \eqdot (\mu_1, \mu_2, \mu_u) \in \Cal P = [6.0, 20.0]  \times [1.0, 3.0]\times (0,1)$. 
The fixed configuration $\Omega_o$ refers to the spatial domain $\Omega$ of Figure \ref{fig:dominio_1} for  a reference parameter $\mu_u^o = 0.5$: we choose the middle point as the most representative case. 

In the Geo-R framework, a new mesh was used, thus, in this case, the high-fidelity fidelity dimension is $2N_h = 10356$. The POD is performed in the reference framework over a uniformly distributed training set of $100$ samples. Indeed, the problem is simpler and less snapshots can be explored to reach good results. Defining the vector of coordinates in the reference domain as $\mathbf x_o = [x_1^o, x_2^o]$, the affine map $T_{\mu_u}: \Omega_o \rightarrow \Omega$ is explicitly defined as follows:
\begin{equation}
  T_{\mu_u}(\mathbf x_o) = \mathbf x =
  \begin{cases}
    \mathbf x_o                & \text{ in } \Omega_1,                        \\
    \begin{bmatrix}(1 + 2(\mu_u - 0.5))(x^o_1 - 1) + 1 \\
      x_2^o
    \end{bmatrix}  & \text{ in } \Omega_i^l, i=2,3,4,\vspace{2mm} \\
    \begin{bmatrix}(1 - 2(\mu_u - 0.5))(x^o_1 - 2) + 2 \\
      x_2^o
    \end{bmatrix} & \text{ in } \Omega_i^r, i=2,3,4.
  \end{cases}
\end{equation}
It is clear that the map $T_{\mu_u}$ verifies the assumptions of Section~\ref{sec:geo-r} and thus the optimization problem can be solved in this new reference domain with numerous benefits.\\
Indeed, solving the problem in $\Omega_o$ allows the state eigenvalues to decay faster than those of the standard POD, as one can see from the right plot of Figure~\ref{fig:error_comp}. Consequently, the average relative error is capable to reach values around $10^{-4}$ for both the variables with only $N=30$ , as depicted in the left plot of Figure~\ref{fig:error_comp}. This result outperforms the standard POD. With respect to L-POD, the state variable is better recovered, while L-POD better reconstructs the adjoint variable up to $N=15$: above that value, the two approaches are totally comparable.

Concluding:
\begin{itemize}
  \item[\small{$\circ$}] standard POD combined with hyper-reduction techniques is an accurate, yet time consuming, way to solve vbOCP($\bmu$)s. Indeed,  a critical amount of basis functions should be used to recover a good representation of the high-fidelity solution. However, no issues are encountered for the adjoint variable.
  \item[\small{$\circ$}] The L-POD helps in reaching accurate results in the online stage. The strength of this strategy is its versatility. It can be applied to every problem one is dealing with. However, it still relies on hyper-reduction techniques as DEIM.
  \item[\small{$\circ$}] A valid option is to recast the problem into an affine decomposed system through Geo-R. Indeed, it allows one to reach more accurate results in a smaller amount of time. The drawback is related to the choice of the reference domain and of the transformation $T_{\mu_u}$. In our case, the choice is made \emph{a posteriori} once observing the physical behaviour of the optimal solutions. We stress that other maps and transformations can be employed. However, this technique can be hardly applied to more complicated geometries.
\end{itemize}

\begin{remark}
  We remark that, in our parametric setting, the Geo-R strategy is morally equivalent to the Arbitrary Lagrangian Eulerian formulation proposed in \cite{torlo2020model}. The reference problem aligns the wave-like phenomenon in the center of the domain and builds a low-dimensional framework eliminating the nonlinear features of the solution manifold. For the sake of completeness, we want to stress that we also tried the \emph{shifted-POD} approach proposed in \cite{nonino2019overcoming}. Despite the promising results they obtain for fluid-structure interaction problems, in our specific context, the performances of the strategy in terms of accuracy were disappointing. {We noticed large errors when we studied limit cases, i.e., $\mu_u \rightarrow 0$ and $\mu_u \rightarrow 1$. This phenomenon was possibly related to interpolation errors and to the naive shift map we used: it aligned the peak at $\mu_u = 0.5$ and extends the solution with the same values along the shifted side. This choice was natural. Nevertheless, it creates fictitious information on the behaviour of the solution and might compromise the accuracy of the ROM solver.}
For this reason we decide not to show the results and we strongly believe that a deeper analysis of such an approach is needed for vbOCP($\boldsymbol \mu$). However, this is beyond the scope of this contribution.
\end{remark}
We stress that the Geo-R allows an affine decomposition of  the system and thus and efficient offline-online decoupling of the process.  The \emph{speed-up index}, i.e.\ the number of reduced problem one can perform in the time of high-fidelity solution is around $92$, while for the standard POD with $N=80$, after the DEIM approximation, the speed-up index is 6. The Geo-R approach is more convenient also with respect to the L-POD: indeed with the same value of $N=30$, L-POD features a speed-up of 17 (the DEIM approach in the online phase affects on the computational performances). Moreover, the advantages of using Geo-R approach for simple geometries is also highlighted from the offline computational costs, since no DEIM is needed in the recast framework. We underline that the small offline cost is also due to the fact that only $100$ snapshots are employed in the GEO-R context. We refer the reader to Table~\ref{tab:case1} for a recap on the  computational times. We stress that, as expected, even if L-POD has larger offline costs with respect to standard POD, they are still acceptable and lead to more accurate and faster results in the online phase.

\begin{table}[]
  \centering
  \caption{Comparison of POD, L-POD and Geo-R in terms of offline computational costs and speed-up.}
  \label{tab:case1}
  {%
    \begin{tabular}{|l|l|l|}
      \hline
      Strategy & Offline Costs (time in seconds)        & Speed-up \\ \hline
      POD      & 935.94s             (DEIM $\sim$ 668s) & 6        \\ \hline
      L-POD    & 1011.49s    (DEIM $\sim$ 668s)         & 17       \\ \hline
      Geo-R    & 68.35s        (NO DEIM needed)         & 92       \\ \hline
    \end{tabular}%
  }
\end{table}

\subsection{Test Case 2: POD and L-POD over a complex geometry}
\label{sec:dolphinPOD}
Despite its capabilities, the Geo-R approach cannot be applied in more complicated settings and it is strictly problem dependent.
Therefore, we propose here a case test where no Geo-R can be employed. This test is of interest for real-life applications, since the vbOCP($\bmu$)s can be related to the simulation of systems characterized by complex geometries, as we said in the Introduction.
The governing PDE($\boldsymbol \mu$) is always equation \eqref{eq:test_1}. 

\begin{figure}[H]
  \centering
  \includegraphics[width=0.4\textwidth]{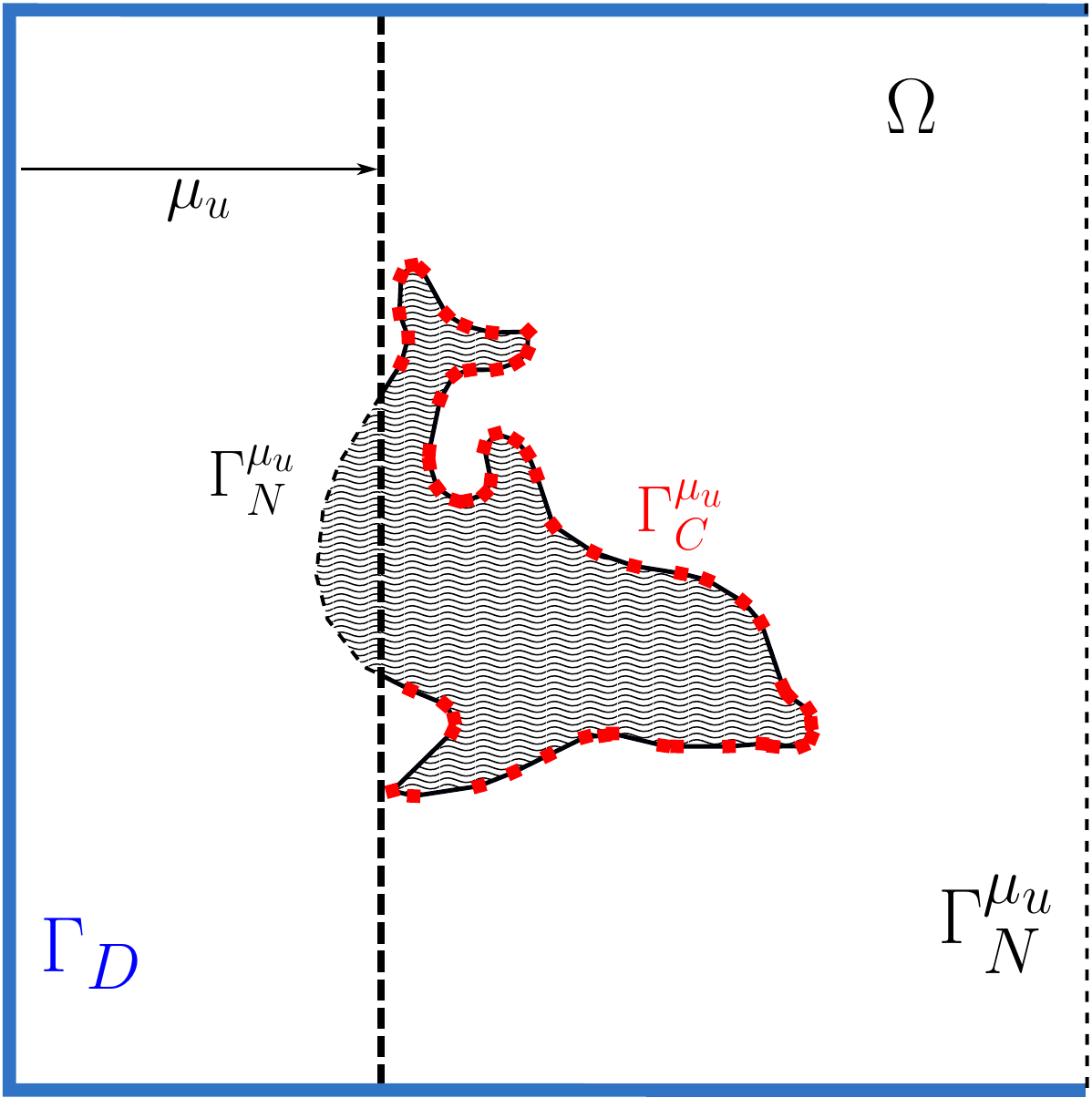}
  \caption{(Test 2). Domain $\Omega$. \textit{Observation domain:} $\Omega_{\text{obs}} = \Omega$, \textit{Control domain:} $\Gamma_C^{\mu_u}$ (red dots). \textit{Blue line:} Dirichlet boundary conditions. \textit{Thin black dashed line:} Neumann boundary conditions. We represent the domain for $\mu_u = 0.33$.}
  \label{fig:dominio_dolph}
\end{figure}

The spatial domain is represented in Figure~\ref{fig:dominio_dolph}. 
We define $\Gamma_{ib}$ as the boundary of the \emph{dolphin hole} in a unit square geometry\footnote{We are deeply thankful to the contributors who provided this specific geometry which is open access to the following link \href{https://fenicsproject.org/pub/data/meshes/dolfinfine.xml}{https://fenicsproject.org/pub/data/meshes/dolfin\textunderscore fine.xml}.}.
In this test case, $\Gamma_C^{\mu_u} \subset \Gamma_{ib}$.
A homogeneous Neumann condition is imposed on the right border of the square and Dirichlet boundary conditions are applied on the other external boundary of the unit square domain.
The parameter is $\bmu \eqdot (\mu_1, \mu_2, \mu_u) \in
  \Cal P = [6.0, 20.0]  \times [0.5, 3.0]\times (0.27, 0.74)$. The parameter $\mu_1$ represents the P\'eclet number and $\mu_2$ describes the desired constant state. This time, the observation is taken all over the domain. The parameter $\mu_u$ ranges from the minimum and the maximum $x_1$ values of $\Gamma_{ib}$.
  For this test case, $\Gamma_C^{\mu_u} = ((\mu_u, 0.74) \times [0,1]) \cap \Gamma_{ib}$ and $\Gamma_N^{\mu_u} = (\Gamma_{ib} \setminus \Gamma_C^{\mu_u}) \cup (\{1\} \times [0, 1])$.
  We set the penalization parameter as the previous experiment: $\alpha = 0.07$.

\begin{figure}
  \centering
  \includegraphics[width=0.48\textwidth]{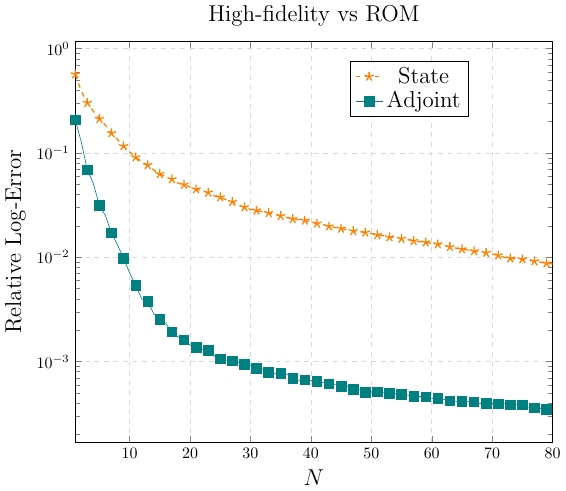}\
  \includegraphics[width=0.48\textwidth]{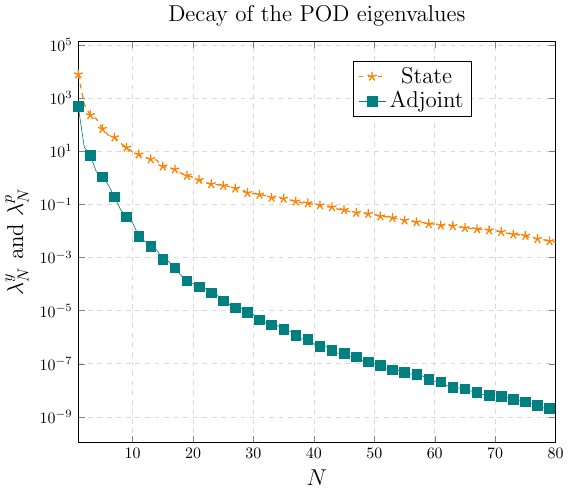}
  \caption{(Test 2: POD). {\emph{Left}. Averaged relative log-error for the two variables. \emph{Right}. Decay of the eigenvalues for the two variables}.}
  \label{fig:dolph_error_lambda300POD}
\end{figure}

\begin{figure}[H]
  \centering
  \includegraphics[width=0.3\textwidth]{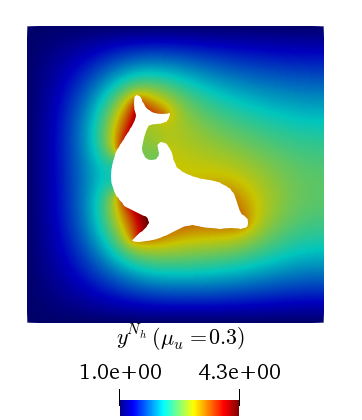}
  \includegraphics[width=0.3\textwidth]{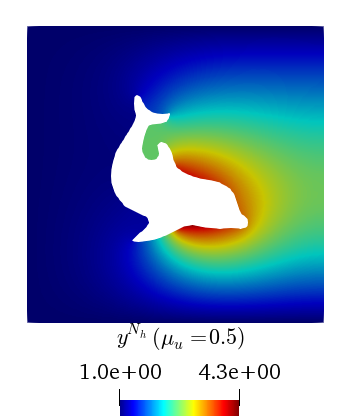}
  \includegraphics[width=0.3\textwidth]{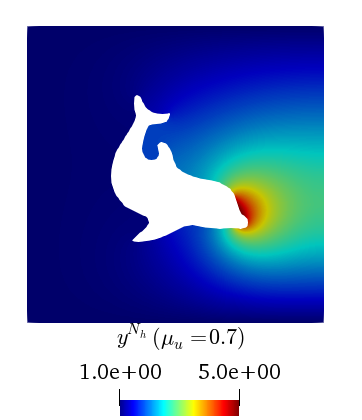}\\
  \includegraphics[width=0.3\textwidth]{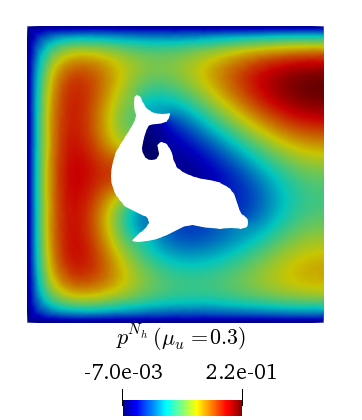}
  \includegraphics[width=0.3\textwidth]{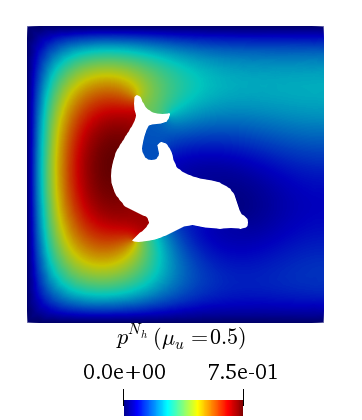}
  \includegraphics[width=0.3\textwidth]{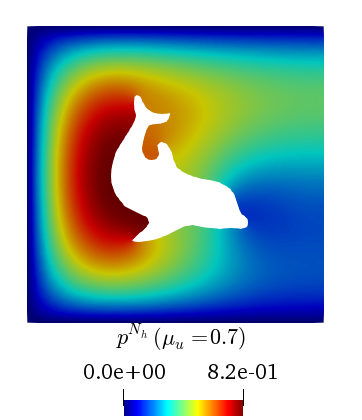}\\
  \caption{(Test 2: POD).  Some high-fidelity solutions for fixed physical parameters ($\mu_1 = 12, \mu_2 = 2.5$) and varying $\mu_u \in \{ 0.3, 0.5, 0.7\}$, from left to right. The state solutions and the adjoint ones are represented on the top and on the bottom of the Figure, respectively.}
  \label{fig:dolph_snap}
\end{figure}

Let us focus on the POD performances. In order to recover the affine decomposition we exploited a DEIM approximation over $350$ uniformly distributed parameters to reach a tolerance of $10^{-5}$ for the eigenvalue decay (i.e.\ 170 DEIM basis). The POD algorithm is ran over another $300$ uniform distributed snapshots. The high-fidelity system is tackled through a $\mathbb{P}^1-\mathbb{P}^1$ approximation of dimension $2N_h = 5736$.
Figure~\ref{fig:dolph_error_lambda300POD} (left plot) shows the relative errors $E_y$ and $E_p$ averaged over a testing set of 150 uniformly distributed parameters in $\mathcal P$.
The complexity of the problem is visible from the plot. Indeed, $80$ basis functions are not enough to reach an accurate state representation, with $E_y \sim 10^{-2}$. While $E_p \sim 10^{-3}$ already for $N=30$. The claim is confirmed by the decay of the eigenvalues represented in the right plot of Figure~\ref{fig:dolph_error_lambda300POD}, which is faster for the adjoint variable. \\
The main features of the problem can be observed in the plots of Figure~\ref{fig:dolph_snap}: we have multiple ``peaks" around the dolphin for small values of $\mu_u$ and they transfer from left to right for the state variable, while the adjoint variable features a very complex behaviour. In this case, there is no symmetry and the geometrical properties are too complex to perform Geo-R. However, we employ L-POD also in this case. As in the previous test case, we built a global DEIM approximation to recover the affinity assumption. As usual, we used 350 uniformly distributed parameters with the usual tolerance of $10^{-5}$. The basis number of the DEIM approximation is $168$ (different from the standard POD since the 350 parameters for the hyper-reduction are different from the ones of the standard POD). We run algorithm \ref{al:01} with $M = 10$, $\tau = 10^{-3}$ and $N=30$: the interval $I_{\mu_u}$ is divided in $J=7$ sub-intervals as
$$
  I_{\mu_u} = \bigcup_{i=1}^J I_{\mu_u}^i,$$
where
$
  I_{\mu_u}^1 = (0.27, 0.32875],  I_{\mu_u}^2 = (0.32875, 0.3875],   I_{\mu_u}^3 = (0.3875, 0.44635] ,
  I_{\mu_u}^4 =(0.44635, 0.505], $
$I_{\mu_u}^5 =  (0.505, 0.6225],   I_{\mu_u}^6 =(0.6225, 0.68125] \text{ and } I_{\mu_u}^7= (0.68125,0.74)
$.

\begin{figure}[H]
  \centering
  \includegraphics[width=0.48\textwidth]{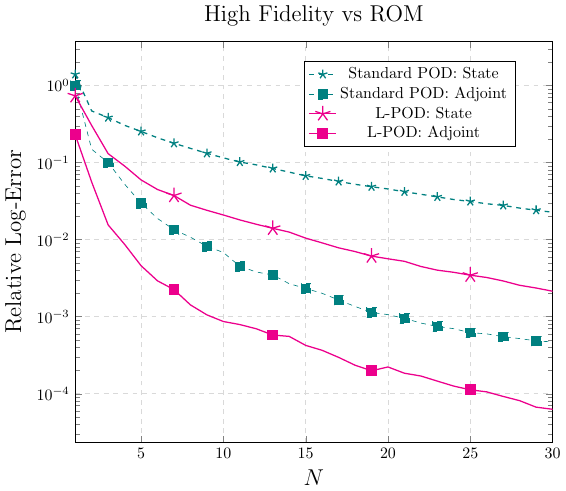}\
  \includegraphics[width=0.48\textwidth]{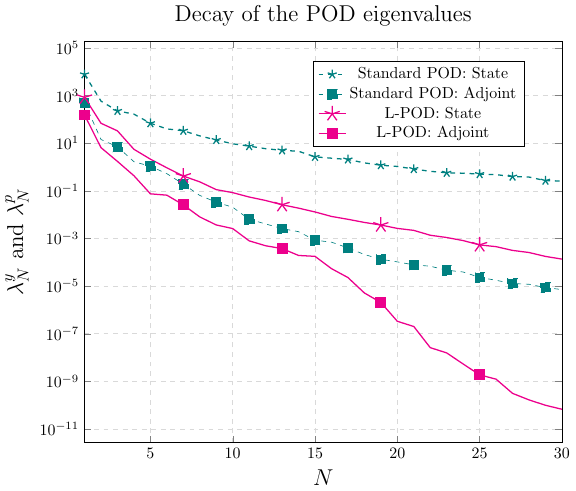}
  \caption{(Test 2: approach comparison). {\emph{Left}. Averaged relative log-error for the two variables for all the approaches. \emph{Right}. Decay of the eigenvalues for the two variables for all the approaches}.}
  \label{fig:dolph_error_lambda300}
\end{figure}

Thanks to the local approach, if we look at the average eigenvalues in the right plot of Figure~\ref{fig:dolph_error_lambda300}, we can observe a better behaviour for both the variables. For this reason, the averaged errors on 150 uniformed distributed parameters is lower with respect to standard POD. With L-POD we reach $E_y \sim 2 \cdot 10^{-3}$ and $E_p \sim 7 \cdot 10^{-5}$, as depicted in the left plot of Figure~\ref{fig:dolph_error_lambda300}. We stress that neither $N=80$ suffices to reach these values of the relative errors. With a comparable offline phase, the L-POD is also convenient in terms of online costs, with a speed-up around 9. The standard POD features a speed-up of only $2$ using $N=80$. We stress that with comparable online costs, i.e.\ for $N=30$, L-POD outperforms standard POD in accuracy.
We recap the computational times in Table~\ref{tab:case2}.

\begin{table}[]
  \centering
  \caption{Comparison of POD and L-POD in terms of offline computational costs and speed-up.}
  \label{tab:case2}
  {%
    \begin{tabular}{|l|l|l|}
      \hline
      Strategy & Offline Costs (time in seconds)     & Speed-up \\ \hline
      POD      & 1844.84    (DEIM $\sim$ 1247s)      & 2        \\ \hline
      L-POD    & 1941.68s       (DEIM $\sim$  1247s) & 9        \\ \hline
    \end{tabular}%
  }
\end{table}

\section{Conclusions}
\label{sec:conclusions}
In this work we deal with ROM techniques to deal with vbOCP($\boldsymbol \mu$)s, i.e.\ boundary optimal control problems where the geometrical action of the control is related to a parameter. Namely, a specific parameter changes the portion of the Neumann boundary where the control plays a role. To the best of our knowledge, it is the first time that this kind of parametric optimal control problem is investigated.

The problem turned out to be difficult to be reduced, since features very complex parametric behaviours. 
Moreover, the vbOCP($\boldsymbol \mu$)s does not verify the affine assumption. We recover it by means of a DEIM algorithm when needed. Guided by the numerical results, we propose tailored approaches inspired by classical methods used in wave-like phenomena and we compared them to standard POD: the Geo-R and the L-POD.
The advantages of such strategies have been tested  on two numerical experiments: a vbOCP($\boldsymbol \mu$) over a simple geometry (where the Geo-R can be applied) and a vbOCP($\boldsymbol \mu$) over a complex one (where only L-POD can be performed and compared to standard POD). The proposed tailored strategies outperform standard POD.

We observe that the reduction of vbOCP($\bmu$)s is a though task due to the transportation of ``peaks" arising for the control action. The diversity of the behaviour of snapshots does not allow a good representation of the physical phenomenon with an acceptable number of basis functions with standard POD technique.
In simple cases, one can rely on Geo-R approximation to tackle a simpler problem, where the control boundary is fixed and thus good performances with a smaller amount of basis are guaranteed. Geo-R guarantees good speed-up results, due to the affine structure of the reference problem. However the choice of the reference domain and the map $T_{\mu_u}$ is strictly problem dependent and made \emph{a posteriori} once observed the problem at hand. 
To overcome this limitation, L-POD can be employed. Indeed, L-POD gives more accurate results with respect to standard POD not paying much effort in the offline computation phase.
This approach is very general and can be applied to complicated problems characterized by high complex geometries. We stress that L-POD still needs hyper-reduction techniques, as standard POD. However, L-POD is beneficial since it recovers the solution variables more accurately with respect to standard POD for the same values of $N$ (i.e.\ for online comparable computational costs).

This contribution is a first step towards many further developments. A possible advance is represented by the analysis of error certification for this specific optimal control formulation. {Indeed, the problem complies with the error analysis proposed in \cite{negri2013reduced} and \cite{StrazzulloRB}, however, it suffers large estimator values for the limit case scenarios. Thus, an investigation of a tailored certification can be important for this particular application.} \\
{A further interesting advancement would be the analysis of 3D complex geometries. We are aware that the adaptive L-POD strategy (even if generalizable to 3D structures) might suffer higher geometrical complexity resulting in many reduced subspaces and many basis functions to deal with. We are confident that intrusive ROM strategies enhanced by cutting-edge techniques (based on machine learning, for example) might be helpful in this context.}\\
Nevertheless, we believe that the present contribution paves the way to some interesting applications in interdisciplinary fields as geophysics, faults and fractures analysis and energy engineering where the need of fast and reliable simulations for forecast intents is daily increasing.

\section*{Acknowledgements}
The computations in this work have been performed with RBniCS \cite{rbnics} library, developed at SISSA mathLab, which is an implementation in FEniCS \cite{fenics} of several reduced order modelling techniques; we acknowledge developers and contributors to both libraries. Computational resources were partially provided by HPC@POLITO, a project of Academic Computing within the Department of Control and Computer Engineering at the Politecnico di Torino (\href{http://hpc.polito.it}{http://hpc.polito.it}). This work is partially supported by the INdAM-GNCS project ``Metodi numerici per lo studio di strutture geometriche parametriche complesse" (CUP\_E53C22001930001) and by the MIUR project ``Dipartimenti di Eccellenza 2018-2022'' (CUP E11G18000350001).

\bibliographystyle{abbrv}
\bibliography{maria.bib}

\end{document}